\documentclass[11pt]{amsart}

\usepackage {enumerate}
\usepackage{setspace}
\usepackage{caption}
\usepackage{mathrsfs}
\usepackage{braket}
\usepackage{slashed}

\usepackage{hyperref}

\newtheorem{theorem}{Theorem}[section]
\newtheorem{lemma}[theorem]{Lemma}
\newtheorem{corollary}[theorem]{Corollary}

\theoremstyle {definition}
\newtheorem{remark}[theorem]{Remark}

\usepackage{tikz}
\usetikzlibrary{decorations.pathreplacing}

\usepackage[margin=1in]{geometry}

\allowdisplaybreaks
  
\DeclareMathOperator{\area}{area} 
\DeclareMathOperator{\vol}{vol}
\newcommand{\R}{\mathbb{R}}
\newcommand{\RR}{\mathbb{R}}
\DeclareMathOperator{\Ric}{Ric}

\DeclareMathOperator{\proj}{proj}
\DeclareMathOperator{\tr}{tr}
\DeclareMathOperator{\Div}{div}

\newcommand{\sF}{\mathscr{F}}

\newcommand{\bangle}[1]{\left\langle #1 \right\rangle}

\title[Far-outlying CMC surfaces] 
{On far-outlying CMC spheres in asymptotically flat Riemannian $3$-manifolds}

\author{Otis Chodosh}
\address{Princeton University, Department of Mathematics, Fine Hall, Washington Road, Princeton, NJ, 08544, United States}
\email{ochodosh@princeton.edu}

\author{Michael Eichmair}
\address{Faculty of Mathematics, University of Vienna, Oskar-Morgenstern-Platz 1, 1090 Vienna, Austria}
\email {michael.eichmair@univie.ac.at}

\begin{document}

\begin{abstract} We extend the Lyapunov-Schmidt analysis of outlying stable CMC spheres in the work of S. Brendle and the second-named author \cite{Brendle-Eichmair:2014} to the ``far-off-center" regime and to include general Schwarzschild asymptotics. We obtain sharp existence and non-existence results for large stable CMC spheres that depend very delicately on the behavior of scalar curvature at infinity. 

\end{abstract}

\maketitle

\date{}


\section {Introduction}

We complement our recent work \cite{angstnomore} on the characterization of the leaves of the canonical foliation as the \emph{unique large closed embedded stable constant mean curvature surfaces} in strongly asymptotically flat Riemannian $3$-manifolds. More precisely, we extend here the Lyapunov-Schmidt analysis of outlying stable constant mean curvature spheres that developed by S. Brendle and the second-named author in \cite{Brendle-Eichmair:2014} to also include the far-off-center regime and general Schwarzschild asymptotics. \\

We begin by introducing some standard notation. \\

Throughout this paper, we consider complete Riemannian $3$-manifolds $(M, g)$ so there are both a compact set $K \subset M$ and a diffeomorphism
\begin{align*}
M \setminus K \cong \{x \in \R^3 : |x| > 1/2\}
\end{align*}
such that in this \emph{chart at infinity}, for some $q > 1/2$ and non-negative integer $k$, 
\begin {align} \label{eqn:AF}
g_{ij} = \delta_{ij} + \tau_{ij}
\end {align}
where 
\[
\partial_I \tau_{ij} = O (|x|^{-q - |I|})
\]
for all multi-indices $I$ of length $|I| \leq k$. Moreover, we require that the boundary $\partial M$ of $M$, if non-empty,  is a minimal surface, and that the components of $\partial M$ are the only connected closed minimal surfaces in $(M, g)$. We say that $(M, g)$ is $C^k$-\emph{asymptotically flat} of rate $q$. 

It is convenient to denote, for $r > 1$, by $S_r$ the surface in $M$ corresponding to the centered coordinate sphere $S_r(0) = \{ x \in \R^3 : |x| = r\}$, and by $B_r$ the bounded open region enclosed by $S_r$ and $\partial M$. Given $A \subset M$, we let
\[
r_0 (A) := \sup \{ r > 1 : B_r \subset A\}.
\]

A particularly important example of an asymptotically flat Riemannian $3$-manifold is Schwarzschild initial data
\[
M = \{x \in \R^3 : |x| \geq m/2\} \qquad \text{ and } \qquad g =  \Big( 1 + \frac{m}{2 |x|} \Big)^4 \sum_{i=1}^3 dx^i \otimes dx^i
\]
where $m > 0$ is the \emph{mass} parameter. 

We say that $(M, g)$ as above is $C^k$-asymptotic to Schwarzschild of mass $m>0$, if, instead of \eqref{eqn:AF}, we have 
\begin{align} \label{eqn:AS}
g_{ij} = \Big(1 + \frac{m}{2|x|}\Big)^{4}\delta_{ij} + \sigma_{ij}
\end{align}
where
\[
\partial_I \sigma_{ij} = O (|x|^{- 2 - |I|})
\]
for all multi-indices $I$ of length $|I| \leq k$. 

The contributions in this paper combined with the key result in \cite{angstnomore} lead to the following theorem.   

\begin{theorem} [\cite{angstnomore}] \label{theo:ANM-glob-v1}
Let $(M,g)$ be a complete Riemannian $3$-manifold that is $C^6$-asymptotic to Schwarzschild of mass $m > 0$ and whose scalar curvature vanishes. 
Every connected closed embedded stable constant mean curvature surface with sufficiently large area is a leaf of the canonical foliation. 
\end{theorem}

The canonical foliation $\{\Sigma_H\}_{0 < H < H_0}$ of $M \setminus K$ through stable constant mean curvature spheres $\Sigma_H$ (with respective mean curvature $H$) was discovered by G.\ Huisken and S.-T.\ Yau \cite{Huisken-Yau:1996}. They show that, for every $s \in (1/2, 1]$, there is $H_s \in (0, H_0)$ such that $\Sigma_H$ for $H \in (0, H_s)$ is the only stable constant mean curvature sphere of mean curvature $H$ in $(M, g)$ that encloses the ball  $\{ x \in \R^3 : |x| < H^{-s} \}$ in the chart at infinity. This characterization of the leaves was later refined by J.\ Qing and G.\ Tian \cite{Qing-Tian:2007}: Upon enlarging $K$ and shrinking $H_0>0$ accordingly, if necessary, each $\Sigma_H$ of the canonical foliation $\{\Sigma_H\}_{0 < H < H_0}$ is the unique stable constant mean curvature sphere of mean curvature $H$ in $(M, g)$ that encloses $K$. In joint work with A. Carlotto \cite{mineffectivePMT}, inspired by earlier work of J. Metzger and the second-named author \cite{stablePMT}, we have extended this characterization further under the additional assumption that the scalar curvature of $(M, g)$ is non-negative in the following way: Choose a point $p \in M$. Every connected stable constant mean curvature sphere $\Sigma \subset M$ that encloses $p$ and whose area is sufficiently large is a leaf of the canonical foliation. Thus, to prove an unconditional uniqueness result along the lines of Theorem \ref{theo:ANM-glob-v1}, it remains to understand large stable constant mean curvature spheres that are \emph{outlying} in the sense that the region they enclose is disjoint and ---  in view of the results in \cite{mineffectivePMT} --- far from $K$. The center of mass flux integrals used in \cite{Huisken-Yau:1996,Qing-Tian:2007} as a centering device vanish in this case regime; new ideas are needed. S. Brendle and the second-named author \cite{Brendle-Eichmair:2014} have discovered a subtle relationship between scalar curvature and outlying stable constant mean curvature spheres. They give examples of divergent sequences $\{\Sigma_k\}_{k = 1}^\infty$ of outlying stable constant mean curvature spheres in $(M, g)$ asymptotic to Schwarzschild with $m>0$, which is the setting of \cite{Huisken-Yau:1996, Qing-Tian:2007}. In fact, $\Sigma_k$ is a perturbation of the coordinate sphere
\[
S_{\lambda_k} (\lambda_k \xi) = \{ |x - \lambda_k \xi| = \lambda_k : x \in \R^3\}
\]
in the chart at infinity, where $\xi \in \R^3$ is such that $|\xi| > 1$ and $\lambda_k \to \infty$.  On the other hand, they show that no such sequences can exist in $(M, g)$ if the scalar curvature is non-negative, provided a further technical assumption on the expansion of the metric in the chart at infinity holds.

\begin{theorem}[S. Brendle and M. Eichmair \cite{Brendle-Eichmair:2014}] \label{thm:BE-non-exist-outlying}
Let $(M,g)$ be a complete Riemannian $3$-manifold that is $C^4$-asymptotic to Schwarzschild with mass $m > 0$, where, in addition to \eqref{eqn:AS}, we also ask that 
\begin{align} \label{eqn:toporderhom}
g_{ij} = \Big( 1 + \frac{m}{2|x|} \Big)^4 \delta_{ij} + T_{ij} + o(|x|^{-2}) \qquad \text{ as } \qquad |x| \to \infty
\end{align}
with corresponding estimates for all partial derivatives of order $\leq 4$, and where $T_{ij}$ is homogeneous of degree $-2$. There does not exist a sequence of outlying stable constant mean curvature surfaces $\{ \Sigma_{k} \subset M \}_{k=1}^\infty$ whose inner radius $r_{0}(\Sigma_{k})$ and mean curvature $H(\Sigma_{k})$ satisfy
\begin{align} \label{conclusionBE}
r_{0}(\Sigma_{k}) \to \infty\qquad \text{and}\qquad r_{0}(\Sigma_{k})H(\Sigma_{k}) \to \eta > 0.
\end{align}
\end{theorem}

In our recent work  \cite{angstnomore}, we show that when $(M, g)$ is asymptotic to Schwarzschild with mass $m > 0$ and if the scalar curvature is non-negative, there are no sequences of embedded stable constant mean curvature spheres $\{\Sigma_k\}_{k=1}^\infty$ in $(M, g)$ with 
\[
r_{0}(\Sigma_{k}) \to \infty \qquad \text{ and } \qquad r_{0}(\Sigma_{k})H(\Sigma_{k}) \to 0.
\]
Assuming in addition that the metric has the form in Theorem \ref{thm:BE-non-exist-outlying}, this leaves only the case of 
\[
r_0(\Sigma_k) \to \infty, \qquad \qquad \area_g (\Sigma_k) \to \infty, \qquad \qquad r_{0}(\Sigma_{k})H(\Sigma_{k}) \to \infty.  
\]
To rule out this scenario, we  revisit the Lyapunov--Schmidt reduction in \cite{Brendle-Eichmair:2014}. 

Our other main goal here is to investigate whether top-order homogeneity in the expansion of the metric \eqref{eqn:toporderhom} off of Schwarzschild in Theorem \ref{thm:BE-non-exist-outlying} is really necessary. Neither the results  \cite{Huisken-Yau:1996,Qing-Tian:2007,mineffectivePMT} for spheres that are not outlying nor the main result of \cite{angstnomore} require such an assumption. It turns out that Theorem \ref{thm:BE-non-exist-outlying}  is \emph{false} without additional such conditions. 

\begin{theorem}\label{thm:exist-outlying}
There is an asymptotically flat complete Riemannian $3$-manifold $(M,g)$ with non-negative scalar curvature that is smoothly asymptotic to Schwarzschild of mass $m>0$ in the sense that
\[
g_{ij} = \Big(1+\frac{m}{2 |x|}\Big)^{4}\delta_{ij} + \sigma_{ij}
\]
where 
\[
\partial_I \sigma_{ij} = O(|x|^{-2 - |I|})
\]
for all multi-indices $I$, which contains a sequence of outlying stable constant mean curvature spheres $\Sigma_{k} \subset M$  with 
\[
r_{0}(\Sigma_{k})\to\infty \qquad \text{ and } \qquad r_{0}(\Sigma_{k}) H(\Sigma_{k}) \to \eta >0.
\] 
\end{theorem}

It turns out that it is possible to recover a version of Theorem \ref{thm:BE-non-exist-outlying} without demanding homogeneity in the expansion of the metric if instead we impose a mild growth condition on the scalar curvature. 

\begin{theorem}\label{thm:non-exist-outlying}
Let $(M,g)$ be a complete Riemannian $3$-manifold that is $C^{4}$-asymptotic to Schwarzschild in the sense that 
\[
g_{ij} = \Big(1 + \frac{m}{2|x|} \Big)^{4}\delta_{ij} + \sigma_{ij}
\]
where $\partial_I \sigma_{ij} = O (|x|^{- 2 - |I|})$ for all multi-indices $I$ of length $|I| \leq 4$. We also assume that either 
\[
R = o (|x|^{-4}) \qquad \text{ as } \qquad |x| \to \infty
\]
or 
\begin{align}\label{radiallnonincreasing}
x^{i}\partial_{i}  ( |x|^2 R ) \leq o(|x|^{-4}) \qquad \text{ as } \qquad |x| \to \infty.
\end{align}
There does not exist a sequence of outlying stable constant mean curvature surfaces $\Sigma_{k} \subset M$ whose inner radius $r_{0}(\Sigma_{k})$ and mean curvature $H(\Sigma_{k})$ satisfy
\[
r_{0}(\Sigma_{k}) \to \infty\qquad \text{and}\qquad r_{0}(\Sigma_{k})H(\Sigma_{k}) \to \eta >0.
\]
\end{theorem}

Note that \eqref{radiallnonincreasing} holds in either one of the following two cases. 
\begin{enumerate} [(i)]
\item When $R = 0$. This is for example the case when $(M,g)$ is \emph{time symmetric} initial data for a vacuum spacetime.
\item When the metric in the chart at infinity has the special form \eqref{eqn:toporderhom} in Theorem \ref{thm:BE-non-exist-outlying}, then 
\[
R = S + o (|x|^{-4}) \qquad \text{ where } \qquad S = \sum_{i, j = 1}^3 \big( \partial_i \partial_j T_{ij} - \partial_i \partial_i T_{jj} \big). 
\]
Note that $S$ is  homogeneous of degree $-4$. Euler's theorem gives that \eqref{radiallnonincreasing} holds if and only if $R  \geq -o(|x|^{-4})$. As such, Theorem \ref{thm:non-exist-outlying} generalizes Theorem \ref{thm:BE-non-exist-outlying} to the non-homogeneous setting.
\end{enumerate}
It is interesting to compare \eqref{radiallnonincreasing} to condition (H3) in S.\ Brendle's version of Alexandrov's theorem for certain warped products  \cite{Brendle:2013}. We remark that the example constructed in Theorem \ref{thm:exist-outlying} is a warped product. 
We also mention that S. Ma has constructed examples of $(M, g)$ that contain large \emph{unstable} constant mean curvature spheres  \cite{Ma:unstable}. The scalar curvature in these examples is negative in some places; see the discussion preceding the statement of Theorem 1.1 in \cite{Ma:unstable} and the proof of Lemma 4.7 therein.

We now turn to the case of surfaces that are \emph{very far} outlying in the sense that 
\[
r_{0}(\Sigma_{k})H(\Sigma_{k})\to\infty.
\]
These surfaces are not within the scope of the Lyapunov--Schmidt reduction carried out in \cite{Brendle-Eichmair:2014}, where the case \eqref{conclusionBE} is considered. 
The main difficulty in this regime is that the ``Schwarzschild contribution'' to the reduced area functional leveraged in \cite{Brendle-Eichmair:2014} is no longer on the order of $O(1)$, but is instead decaying. As such, it is necessary to obtain rather involved estimates for the reduced functional. To describe our results, we first  recall some terminology from \cite{Brendle-Eichmair:2014} that we will also adopt. A standard application of the implicit function theorem gives that, for $\lambda > 0$ and $\xi \in \R^3$  large, we can find closed surfaces $\Sigma_{(\xi,\lambda)}$ in the chart at infinity so that the following hold:
\begin{itemize}
\item $\Sigma_{(\xi,\lambda)}$ bounds volume $4\pi\lambda^{3}/3$ with respect to the metric $g$.
\item $\Sigma_{(\xi,\lambda)}$ is the Euclidean graph of a function $u_{(\xi,\lambda)}$ on $S_{\lambda}(\lambda \, \xi)$, i.e.
\[
\Sigma_{(\xi, \lambda)} = \{ \lambda \, \xi + y + u_{(\xi, \lambda)} (x) \, y / \lambda : x = \lambda \, \xi + y  \in S_{\lambda} (\lambda \, \xi)\},
\] 
where 
\[
\sup_{S_{\lambda}(\lambda \, \xi)}  |u_{(\xi,\lambda)}| + \lambda \sup_{S_{\lambda}(\lambda \, \xi)}   |\nabla u_{(\xi,\lambda)}| + \lambda^2 \sup_{S_{\lambda}(\lambda \, \xi)}  | \nabla^2 u_{(\xi,\lambda)} |  = O (1 / |\xi|).
\]
\item $u_{(\xi, \lambda)}$ is orthogonal to the first spherical harmonics on $S_{\lambda} ( \lambda \, \xi)$ with respect to the Euclidean metric.
\item The mean curvature of $\Sigma_{(\xi,\lambda)}$ with respect to $g$ viewed as a function on $S_{\lambda}(\lambda \, \xi)$ is the restriction of a linear function. 
\end{itemize}
Given a sequence of connected closed stable constant mean curvature surfaces $\{ \Sigma_{k} \}_{k=1}^\infty$ with $r_{0}(\Sigma_{k})\to\infty$ and $r_{0}(\Sigma_{k})H(\Sigma_{k})\to\infty$, the same argument as in \cite[p.\ 676]{Brendle-Eichmair:2014} shows that $\Sigma_k = \Sigma_{(\xi_k, \lambda_k)}$ for appropriate $\lambda_k > 0$ and $\xi_k \in \R^3$ when $k$ is sufficiently large. Note that $\lambda_k > 0$ and $\xi_k \in \R^3$ are both large in this case.  Whether $(M, g)$ admits such sequences can now be decided using the following result. 

\begin{theorem}\label{theo:far-off-center-LS}
Let $(M,g)$ be a complete Riemannian $3$-manifold that is $C^{5+\ell}$-asymptotic to Schwarzschild with mass $m=2$, where $\ell \geq 0$ is an integer. Let $\lambda > 0$ and $\xi \in \R^3$ be large. We have\footnote{We may compute the Laplacian of scalar curvature either with respect to $g$ or with respect to the Euclidean background metric in the chart at infinity. The difference may be absorbed into the error terms of the expansion.}  
\begin{align} \label{lsreduction}
\area_ g(\Sigma_{(\xi,\lambda)}) = 4\pi \lambda^{2} - \frac{2\pi}{15} \lambda^{4} R(\lambda \, \xi) - \frac{\pi}{105} \lambda^{6} (\Delta R)(\lambda \, \xi) - \frac{8\pi}{35}|\xi|^{-6} + O(\lambda^{-1}|\xi|^{-6}) + O(|\xi|^{-7})
\end{align}
where $R$ is the scalar curvature  of $(M, g)$. This expansion can be differentiated $\ell$ times with respect to $\xi$.  
\end{theorem}

As in \cite{Brendle-Eichmair:2014}, we use that for $\ell \geq 1$, the map 
\[
\xi \mapsto \area_g (\Sigma_{(\xi,\lambda)})
\]
has a critical point at $\xi$ if, and only if, $\Sigma_{(\xi,\lambda)}$ is a constant mean curvature sphere. If $\ell \geq 2$, then the critical point is stable if, and only if, $\Sigma_{(\xi,\lambda)}$ is a stable constant mean curvature sphere. This immediately leads to the following corollary. 

\begin{corollary}
Let $(M,g)$ be a complete Riemannian $3$-manifold that is $C^{6}$-asymptotic to Schwarzschild in the sense that
\[
g_{ij} = \Big(1+\frac{m}{2|x|}\Big)^{4}\delta_{ij} + \sigma_{ij}
\]
where $\partial_I \sigma_{ij} = O(|x|^{-2 - |I|})$ for all multi-indices $I$ of length $|I| \leq 6$. Assume that the scalar curvature vanishes. There does not exist a sequence of connected closed stable constant mean curvature surfaces $\{ \Sigma_{k} \}_{k = 1}^\infty$ in $(M, g)$ with 
\[
r_{0}(\Sigma_{k})\to\infty,  \qquad \area_g(\Sigma_k) \to \infty, \qquad \text{and} \qquad r_{0}(\Sigma_{k})H(\Sigma_{k})\to\infty.
\]
\end{corollary}

Lyapunov-Schmidt reduction has also been used by e.g. R. Ye \cite{Ye:1991}, S. Nardulli \cite{Nardulli:2009}, and F. Pacard and X. Xu in \cite{PacardXu} to study when small geodesic spheres admit perturbations to constant mean curvature. S. Nardulli \cite{Nardulli:2009} has studied the expansion for small volumes of the isoperimetric profile of a Riemannian manifold. \\  

The analogue of Theorem \ref{thm:non-exist-outlying} in this setting is not so clear-cut. We have the following result. 

\begin{corollary}\label{cor:veryfar-concave}
Let $(M,g)$ be a complete Riemannian $3$-manifold that is $C^{7}$-asymptotic to Schwarzschild in the sense that
\[
g_{ij} = \Big(1+\frac{m}{2|x|}\Big)^{4}\delta_{ij} + \sigma_{ij}
\]
where $\partial_I \sigma_{ij} = O(|x|^{-2 - |I|})$ for all multi-indices $I$ of length $|I| \leq 6$. We also assume that the scalar curvature $R$ of $(M, g)$ is radially convex at infinity in the sense that  
\begin{equation}\label{eq:R-concave-radial}
x^{i}x^{j}\partial_i \partial_j R \geq 0
\end{equation} 
outside of a compact set. 
There does not exist a sequence of connected closed stable constant mean curvature surfaces $\{ \Sigma_{k} \}_{k = 1}^\infty$ in $(M, g)$ with 
\[
r_{0}(\Sigma_{k})\to\infty,  \qquad \area_g(\Sigma_k) \to \infty, \qquad \text{and} \qquad r_{0}(\Sigma_{k})H(\Sigma_{k})\to\infty.
\]
\end{corollary}

It turns out that the hypothesis \eqref{eq:R-concave-radial} is  surprisingly sharp. Comparing with Theorem \ref{thm:BE-non-exist-outlying} or Theorem \ref{thm:non-exist-outlying}, one might be lead to conjecture that it can be weakened to 
\begin{enumerate} [(i)]
\item assuming that $x^{i}x^{j} \partial_i \partial_j R \geq -o(|x|^{-4})$ as $|x| \to \infty$, or
\item assuming that $\sigma_{ij}=T_{ij}+o(|x|^{-2})$ as $|x| \to \infty$ where $T_{ij}$ homogeneous of order $-2$, and that the scalar curvature is non-negative.
\end{enumerate}
The second alternative assumption here implies the first --- by Euler's theorem. 

The following example dashes any hope of such generalizations. 

\begin{theorem}\label{thm:homog-counterexample}
There is an asymptotically flat complete Riemannian $3$-manifold $(M,g)$ with non-negative scalar curvature such that, in the chart at infinity, 
\[
g_{ij}=(1+|x|^{-1})^{4}\delta_{ij} + T_{ij} + o(|x|^{-4}) \qquad \text{ as } \qquad |x| \to \infty 
\]
along with all derivatives, where $T_{ij}$ is homogeneous of degree $-2$, and which contains  outlying stable constant mean curvature spheres $\Sigma_{k} \subset M$ with 
\[
r_{0}(\Sigma_{k}) \to\infty, \qquad \qquad \area_g(\Sigma_k) \to \infty,  \qquad  \text{ and } \qquad r_{0}(\Sigma_{k})H(\Sigma_{k})\to\infty.
\] 
\end{theorem}

Finally, we note that there is by now an impressive body of work on stable constant mean curvature spheres in general asymptotically flat Riemannian $3$-manifolds. We refer the reader to Section 2.1 in \cite{angstnomore} for an overview and references to results in this direction. 


\subsection*{Acknowledgments} 

We thank S. Brendle for helpful conversations. M. Eichmair has been supported by the START-Project Y963-N35 of the Austrian Science Fund. 


\section{Proof of Theorem \ref{thm:non-exist-outlying}}

The proof follows the strategy of \cite{Brendle-Eichmair:2014}, with one important difference: We do not assume here that the deviation of the metric from Schwarzschild is homogeneous of degree $-2$ to top order. Without loss of generality, we may assume that the mass $m$ is equal to $2$. Thus,
\[
g_{ij} =(1+|x|^{-1})^{4}\delta_{ij} + \sigma_{ij}
\]
where 
\[
\partial_{I}\sigma_{ij} = O(|x|^{-2-|I|})
\]
for all multi-indices $I$ of length $|I|\leq 4$. 

Let $\Omega$ be a bounded subset with compact closure in $\RR^{3}\setminus \overline{ B_{1}(0)}$. For $\xi \in \Omega$ and $\lambda >0$ sufficiently large, we may use the implicit function theorem to find surfaces $\Sigma_{(\xi,\lambda)}$ as in Proposition 4 of \cite{Brendle-Eichmair:2014}. Moreover, the surface $\Sigma_{(\xi,\lambda)}$ is a constant mean curvature sphere (respectively, a stable constant mean curvature sphere) if, and only if, $\xi$ is a critical point (respectively, a stable critical point) for the map 
\begin {align*} 
\xi \mapsto \area_{g}(\Sigma_{(\xi,\lambda)}).
\end {align*}
The derivation of Proposition 5 in \cite{Brendle-Eichmair:2014} carries over to give 
\begin{align} \label{BEexpansion}
\area_{g}(\Sigma_{(\xi,\lambda)})  = 4\pi \lambda^{2} + \frac \pi 2  F_{0}(|\xi|) +  F_{\sigma}(\xi,\lambda) + o(1) \qquad \text{ as } \qquad \lambda \to \infty.
\end{align}
The assumption that $\sigma$ is homogeneous is neither needed nor used at this point of \cite{Brendle-Eichmair:2014}. We recall that  
\[
F_{0}(t) = -14 + 16 \, t^{2} \log \frac{t^{2}-1}{t^{2}} + (15 \, t -t^{-1}) \log\frac{t+1}{t-1}
\]
is the contribution from the Schwarzschild background, while 
\begin{align} \label{eqn:reduction}
F_{\sigma}(\xi,\lambda) = \frac 12 \int_{S_{(\xi,\lambda)}}  \tr_{S_{(\xi,\lambda)}} \sigma  - \frac 1 \lambda \int_{B_{(\xi,\lambda)}}  \tr \,  \sigma  
\end{align}
is the contribution from $\sigma$. 

Here and below, unless explicitly noted otherwise, all geometric operations are with respect to the Euclidean background metric in the chart at infinity. 

As in \cite{Brendle-Eichmair:2014}, given $\xi \in \R^3$ and $\lambda > 0$, we will often write 
\[
S_{(\xi, \lambda)} = S_\lambda (\lambda \, \xi) \qquad \text{ and } \qquad B_{(\xi, \lambda)} = B_\lambda (\lambda \, \xi).
\]


\subsection{Radial variation} \label{subsec:radial-1st-var}

The computation of the radial derivative of \eqref{eqn:reduction} in Section 3 of \cite{Brendle-Eichmair:2014} uses the top-order homogeneity of $\sigma$ that is part of their assumption repeatedly. Here, we compute this derivative in the general case, employing several integration by parts to derive a geometric expression involving the scalar curvature \emph{on the nose}. \\
\begin{align*}
(\nabla_{\xi}F_{\sigma} )(\xi,\lambda) & = \frac{d}{ds}\Big|_{s=1}F_{\sigma}(s \, \xi,\lambda) \\
& = \frac \lambda 2  \int_{S_{(\xi,\lambda)}}  \tr_{S_{(\xi,\lambda)}} \nabla_{\xi}\sigma   - \int_{S_{(\xi,\lambda)}} ( \tr \,  \sigma) \bangle { \xi, \nu}  
\intertext{We write $\xi=\xi^\top + \bangle{ \xi, \nu} \nu$.}
& = \frac \lambda 2  \int_{S_{(\xi,\lambda)}} ( \tr_{S_{(\xi,\lambda)}} \nabla_{\nu}\sigma) \bangle{\xi,\nu}   - \int_{S_{(\xi,\lambda)}} ( \tr \,  \sigma) \bangle{\xi,\nu}  \\
& \qquad + \frac \lambda 2  \int_{S_{(\xi,\lambda)}} ( \tr_{S_{(\xi,\lambda)}} \nabla_{\xi^\top}\sigma) \\
& = \frac \lambda 2  \int_{S_{(\xi,\lambda)}} (\tr_{S_{(\xi,\lambda)}} \nabla_{\nu}\sigma) \bangle{\xi, \nu}   - \int_{S_{(\xi,\lambda)}} ( \tr \,  \sigma) \bangle{\xi, \nu} \\
& \qquad + \frac \lambda 2  \int_{S_{(\xi,\lambda)}} (\nabla_{\xi^\top}  \tr \,   \sigma - (\nabla_{\xi^\top}\sigma ) (\nu,\nu) \\
& = \frac \lambda 2  \int_{S_{(\xi,\lambda)}} ( \tr_{S_{(\xi,\lambda)}} \nabla_{\nu}\sigma) \bangle{\xi, \nu}   - \int_{S_{(\xi,\lambda)}} ( \tr \,  \sigma) \bangle{\xi, \nu} \\
& \qquad + \frac \lambda 2  \int_{S_{(\xi,\lambda)}} (\nabla_{\xi^\top} ( \tr \,   \sigma - \sigma(\nu,\nu)) +2 \,  \sigma (\nabla_{\xi^\top} \nu,\nu) \\
& = \frac \lambda 2  \int_{S_{(\xi,\lambda)}} ( \tr_{S_{(\xi,\lambda)}} \nabla_{\nu}\sigma) \bangle{\xi, \nu}   - \int_{S_{(\xi,\lambda)}} ( \tr \,  \sigma) \bangle{\xi, \nu} \\
& \qquad + \frac \lambda 2  \int_{S_{(\xi,\lambda)}} ( ( \tr \,   \sigma - \sigma(\nu,\nu))(-\Div_{S_{(\xi,\lambda)}}\xi^\top) +2 \,  \sigma (\nabla_{\xi^\top}  \nu, \nu) \\
& = \frac \lambda 2  \int_{S_{(\xi,\lambda)}} (  \tr_{S_{(\xi,\lambda)}}  \nabla_{ \nu}\sigma) \bangle{\xi, \nu}   - \int_{S_{(\xi,\lambda)}} (  \tr \,  \sigma) \bangle{\xi, \nu} \\
& \qquad +  \int_{S_{(\xi,\lambda)}}  (  \tr \,   \sigma) \bangle{\xi, \nu} - \sigma( \nu, \nu)  \bangle{\xi, \nu} + \sigma (\xi^\top, \nu) \\
& = \frac \lambda 2  \int_{S_{(\xi,\lambda)}} (  \tr_{S_{(\xi,\lambda)}}  \nabla_{ \nu}\sigma) \bangle{\xi, \nu}  +  \int_{S_{(\xi,\lambda)}} (\sigma (\xi,  \nu) -2 \,  \sigma( \nu, \nu)  \bangle{\xi, \nu} ).
\end{align*}
We define a vector field 
\[
Y = \bangle{\xi, \nu} \sigma( \nu,\, \cdot \, )^{ \sharp}
\]
on $S_{(\xi, \lambda)}$ and compute
\[
  \Div_{S_{(\xi,\lambda)}} Y = \frac 1 \lambda \sigma(\xi, \nu) - \frac 1 \lambda \bangle{\xi, \nu}  \sigma( \nu, \nu) + \bangle{\xi, \nu}   \tr_{S_{(\xi,\lambda)}} ( \nabla_{\, \cdot \, }\sigma)( \nu,\, \cdot \, ) + \frac 1\lambda \bangle{\xi, \nu}  \tr_{S_{(\xi,\lambda)}}\sigma.
\]
The first variation formula gives 
\[
 \frac\lambda 2 \int_{S_{(\xi,\lambda)}} (   \tr_{S_{(\xi,\lambda)}} (  \nabla_{\, \cdot \, }\sigma)( \nu,\, \cdot \, )) \bangle{\xi, \nu} = \frac 1 2 \int_{S_{(\xi,\lambda)}} \big(  3 \, \sigma( \nu, \nu) -  \tr_{S_{(\xi,\lambda)}}\sigma\big)  \bangle { \xi, \nu} - \frac 12 \int_{S_{(\xi,\lambda)}} \sigma(\xi, \nu).
\]
We insert this into the above expression, and continue.
\begin{align*}
\frac{d}{ds}\Big|_{s=1}F(s \, \xi,\lambda) & = \frac \lambda 2  \int_{S_{(\xi,\lambda)}} (  \tr_{S_{(\xi,\lambda)}}  \nabla_{ \nu}\sigma -  \tr_{S_{(\xi,\lambda)}} ( \nabla_{\, \cdot \, }\sigma )( \nu,\, \cdot \, ) ) \bangle{\xi, \nu} \\
& \qquad  - \frac 12  \int_{S_{(\xi,\lambda)}} (   \tr \,  \sigma )  \bangle{\xi, \nu} - \sigma (\xi,  \nu) 
\intertext{We write $  \bangle { \xi, \nu} = - |\xi|^{2} + \lambda^{-1}  \bangle { \xi,X}$ in the first integrand, where $X$ is the position field.}
& =  \frac \lambda 2  \int_{S_{(\xi,\lambda)}} ( \tr_{S_{(\xi,\lambda)}} ( \nabla_{\, \cdot \, }\sigma )( \nu,\, \cdot \, ) -   \tr_{S_{(\xi,\lambda)}}  \nabla_{ \nu}\sigma ) (|\xi|^{2} -\lambda^{-1}  \bangle { \xi,X}) \\
& \qquad  - \frac 12  \int_{S_{(\xi,\lambda)}} ( (   \tr \,  \sigma )  \bangle{\xi, \nu} - \sigma (\xi,  \nu) )  \\
& =  \frac \lambda 2  \int_{S_{(\xi,\lambda)}} (  \tr \, (   \nabla_{\, \cdot \, }\sigma)( \nu,\, \cdot \, ) -   \tr \,    \nabla_{ \nu}\sigma ) (|\xi|^{2} -\lambda^{-1}  \bangle { \xi,X}) \\
& \qquad  - \frac 12  \int_{S_{(\xi,\lambda)}}  (   \tr \,  \sigma )  \bangle{\xi, \nu} - \sigma (\xi,  \nu) 
\intertext{We define a vector field $W =   \Div \sigma -  \nabla  \tr \,  \sigma$.}
& =  \frac 1 2  \int_{S_{(\xi,\lambda)}}   \bangle {\xi, \lambda \, \xi - X} \bangle { W, \nu} \\
& \qquad  - \frac 12  \int_{S_{(\xi,\lambda)}} (    \tr \,  \sigma )  \bangle{\xi, \nu} - \sigma (\xi,  \nu)   \\
& =  \frac 1 2  \int_{B_{(\xi,\lambda)}}   \Div( \bangle {\xi, \lambda \, \xi - X}W) \\
& \qquad  - \frac 12  \int_{S_{(\xi,\lambda)}} ( (   \tr \,  \sigma )  \bangle{\xi, \nu} - \sigma (\xi,  \nu) )  \\
& =  \frac 1 2  \int_{B_{(\xi,\lambda)}} (  \Div \,  W) \bangle {\xi, \lambda \, \xi - X}\\
& \qquad - \frac 1 2  \int_{B_{(\xi,\lambda)}}   \bangle { \xi,W} \\
& \qquad  - \frac 12  \int_{S_{(\xi,\lambda)}}  (   \tr \,  \sigma )  \bangle{\xi, \nu} - \sigma (\xi,  \nu) 
\intertext{Note that $  \bangle { \xi,W} =   \Div (\sigma(\xi,\, \cdot \, ) - ( \tr \,   \sigma) \xi)$. We apply the divergence theorem.}
& =  \frac 1 2  \int_{B_{(\xi,\lambda)}} (  \Div  \, W)  \bangle {\xi, \lambda \, \xi - X} \\
& \qquad - \frac 1 2  \int_{S_{(\xi,\lambda)}} \sigma(\xi, \nu) - (  \tr \,  \sigma) \bangle{\xi, \nu} \\
& \qquad  - \frac 12  \int_{S_{(\xi,\lambda)}}  (  \tr \,  \sigma )  \bangle{\xi, \nu} - \sigma (\xi,  \nu) \\
& =  \frac 1 2  \int_{B_{(\xi,\lambda)}} (  \Div  \, W)  \bangle {\xi, \lambda \, \xi - X}.
\end {align*}
Note that 
\[
 \Div \,  W = R + O(|x|^{-5})
\]
where $R$ is the scalar curvature of $g$. In conclusion, we obtain
\begin {align} \label{eqn:radialderivative}
(\nabla_{\xi}F_{\sigma}) (\xi, \lambda) = \frac{d}{ds}\Big|_{s=1}F(s \, \xi,\lambda) = \frac{1}{2} \int_{B_{(\xi,\lambda)}} \bangle{\xi, \lambda \, \xi - X} R   + o(1) \qquad \text{ as } \lambda \to \infty.
\end{align}
This computation connects the radial derivative of $F_{\sigma}$ with the scalar curvature $R$ of $g$. We emphasize again that our derivation parallels the proof of Proposition 7 in \cite{Brendle-Eichmair:2014}, though we do not assume the top order homogeneity of $\sigma$. 


\subsection{Radial variation in spherical coordinates} 

Assume first that 
\[
R \geq 0 \qquad \text{ and } \qquad x^i \partial_i (|x|^2 R) \leq 0.
\]
For definiteness, we assume that 
\[
\xi=|\xi| \, e_3
\]
where $|\xi| > 1$. In this subsection,  we compute the radial variation 
\[
\int_{B_{\lambda} (\lambda \, \xi) } \bangle {\lambda \, \xi - X, \xi }   R  
\]
in spherical 
\[
(\rho, \phi, \theta) \mapsto (\rho \,  \sin \phi \,  \cos\theta, \rho \,  \sin \phi \, \sin\theta,\rho \, \cos\phi).
\]
on the complement of the $z$-axis. The radial line in direction 
\[
( \sin \phi \cos\theta,  \sin \phi\sin\theta, \cos\phi)
\]
intersects the sphere $B_\lambda (\lambda \, \xi)$ in the $\rho$-interval whose endpoints are the solutions
\[
\rho_{\pm} = \lambda \, |\xi| \,  \big( \cos\phi \pm ( 1/|\xi|^2 - \sin^2 \phi )^{1/2}\big)
\]
of the quadratic equation 
\[
\rho^{2} - 2 \,  \rho \,  \lambda \,  |\xi| \,  \cos\phi + \lambda^{2} \, (|\xi|^{2}-1)  = 0. 
\]
The intersection is non-empty for angles $\phi \in [0, \phi_+]$ where $\phi_+ \in (0, \pi)$ solves 
\[
\sin^{2}\phi_{+} =1 / |\xi|^2.
\]
We then have that 
\begin{align*}
& \int_{B_{\lambda} (\lambda \, \xi)} \bangle {\lambda \, \xi - X, \xi }  R  \\
& \qquad = \int_{0}^{2\pi} \int_{0}^{\phi_{+}} \int_{\rho_{-}}^{\rho_{+}} R(\rho,\phi,\theta) \left( \lambda |\xi|^{2} - \rho  |\xi| \cos\phi \right) \rho^{2} \sin\phi \, d\rho \, d\phi \,  d\theta\\
& \qquad =  |\xi| \int_{0}^{2\pi} \int_{0}^{\phi_{+}} \int_{\rho_{-}}^{\rho_{+}} \rho^2  R(\rho,\phi,\theta) \left(\lambda |\xi| - \rho \cos\phi \right)  \sin\phi \, d\rho \, d\phi  \, d\theta\\
& \qquad \geq  |\xi| \int_{0}^{2\pi} \int_{0}^{\phi_{+}}  ( |\xi| / \cos \phi)^2 R(|\xi|/\cos\phi,\phi,\theta)   \left(  \int_{\rho_{-}}^{\rho_{+}}  \left(\lambda |\xi| - \rho \cos\phi \right) d \rho \right) \sin\phi  \, d\phi \,  d\theta.
\end{align*}
Now, for every $\phi \in (0, \phi_+)$,  
\begin{align*}
\int_{\rho_{-}}^{\rho_{+}} (\lambda|\xi| - \rho \cos \phi) d\rho & = (\rho_+ - \rho_-) \lambda \, |\xi| \, \sin^2 \phi > 0
\end{align*}
so that, in conclusion, 
\[
 \int_{B_{\lambda} (\lambda \, \xi)} \bangle {\lambda \, \xi - X, \xi }  R  \geq 0.
\]
Arguing as in \cite[p.\ 677]{Brendle-Eichmair:2014} shows that $\Sigma_{(\xi,\lambda)}$ cannot be a constant mean curvature sphere.

We now observe that the above arguments go through under the weaker assumption \eqref{radiallnonincreasing}. 
Indeed, using that $R = O (|x|^{-4})$ from asymptotic flatness, we obtain upon integrating inwards from infinity that 
\[
R \geq - o (|x|^{-4}) \qquad \text{ as } \qquad  |x| \to \infty.
\]
Under these assumptions, the preceding computation leads to the estimate 
\[
 \int_{B_{\lambda} (\lambda \, \xi)} \bangle {\lambda \, \xi - X, \xi }  R  \geq  - o(1) \qquad \text{ as } \qquad \lambda \to \infty.
\]
We also mention that \eqref{radiallnonincreasing} is implied by the assumption   
\[
R \geq - o (|x|^{-4}) \qquad \text{ and } \qquad  4 \, R + x^i \partial_i R \leq o (|x|^{-4})   
\]
both as $|x| \to \infty$. In particular, it follows from the assumptions in Theorem \ref{thm:BE-non-exist-outlying}.


\section{Proof of Theorem \ref{thm:exist-outlying}}

Our strategy here parallels the proof of Theorem 1 in \cite{Brendle-Eichmair:2014} in that we construct our metric to have a pulse in its scalar curvature, which in turn forces the reduced area functional $\xi \mapsto \area_g(\Sigma_{(\xi, \lambda)})$ to have stable critical points. Unlike in \cite{Brendle-Eichmair:2014}, our examples are spherically symmetric (which also simplifies the analysis) and, more importantly, they have non-negative scalar curvature. \\

Let $S : (0, \infty) \to (-\infty, 0]$ be a smooth function with 
\[
S^{(\ell)}(r) = O(r^{-4 - \ell}).
\]
We define a smooth function $\varphi : (0, \infty) \to \R$ by 
\begin{align*}
\varphi(r) = \frac{1}{r}\int_{r}^{\infty}(\rho-r) \, \rho \,  S( \rho ) \, d \rho.
\end{align*}
Note that 
\[
\varphi'(r) = - \frac{1}{r^2}\int_{r}^{\infty} \rho^{2} \, S (\rho) \, d \rho
\]
so 
\begin{equation}\label{eq:lap-varphi}
 (r^{2}\varphi')' / r^2 = S(r).
\end{equation}


\begin{lemma}\label{lemm:order-falloff-varphi-counterexample}
We have that  
\[
 \varphi^{(\ell)}(r) = O(r^{-2-\ell}).
\]

\begin{proof}
Because $S(r) = O(r^{-4})$, we see that
\[
\varphi (r) = O (r^{-2}) \qquad \text{ and } \qquad \varphi'(r) = O (r^{-3}).
\]
Using \eqref{eq:lap-varphi}, we find  
\[
\varphi''(r) + 2 \varphi'(r)/r = S(r).
\]
From this, the asserted decay of the higher derivatives can be verified by induction.
\end{proof}
\end{lemma}

On $\R^3 \setminus \{0\}$, we define a conformally flat Riemannian metric 
\[
g = ( 1 + 1/r +\varphi(r)  )^{4}\bar g = (1 + 1/r )^4 \bar g + O (1/r^{2}) 
\]
where $r = |x|$. Note that $g$ is smoothly asymptotic to Schwarzschild with mass $2$. Its scalar curvature is easily computed as  
\[
R =  - 8  ( 1 + 1/r + \varphi(r)  )^{-5} (r^{2} \varphi')' / r^2 =  - 8  (1+O(1/r))  S(r).
\]
In particular, it is non-negative on the complement of a compact set. We now make a particular choice for $S$. Fix  $\chi \in C^{\infty}(\RR)$ that is positive on $(3, 4)$ and suppored in $[3,4]$. Let 
\[
S(r) = - A \sum_{k=0}^{\infty} 10^{-4k}\chi(10^{-k}r)
\]
where $A > 0$ is a large constant that we will fix later. Recall from \eqref{BEexpansion} that 
\[
\area_{g}(\Sigma_{(\xi,\lambda)}) = 4\pi \lambda^{2} + 2\pi F_{0}(|\xi|) + \frac{1}{2\pi} F_{\sigma}(\xi,\lambda) + o(1) \qquad \text{ as } \qquad \lambda \to \infty.
\]
We choose $\xi \in\RR^{3}$ with $2 \leq |\xi| \leq 9$ and $\lambda = 10^{j}$ where $j \geq 1$ is a large integer. Using \eqref{eqn:radialderivative}, we compute the radial derivative as 
\begin{align*}
\frac{d}{ds} \Big|_{s=1} \area_{g}(\Sigma_{(s\xi,\lambda)}) & = 2\pi |\xi| F_{0}'(|\xi|) + \frac{1}{4\pi} \int_{X \in B_\lambda (\lambda \, \xi) }  R (X)\, \bangle {\xi, \lambda \,  \xi - X }   + o(1)\\
& = 2\pi |\xi| F_{0}'(|\xi|) - \frac{2}{\pi} \int_{X \in B_\lambda (\lambda \, \xi)} S(|X|) \bangle {\xi, \lambda \, \xi - X }   + o(1)\\
& = 2\pi |\xi| F_{0}'(|\xi|) + \frac{2A}{\pi} \int_{Y \in B_{1}(\xi)} \chi(|Y|) \bangle{ \xi, \xi - Y}   + o(1) \qquad \text{ as } \qquad \lambda \to \infty.
\end{align*}
When $|\xi| = 2 \sqrt 2$, the integral on the last line is negative. We choose $A > 0$ large so that the sum of the first two terms is negative. When $|\xi| =5$, the second term vanishes while the first term is strictly positive. Thus, for $j \geq 1$ sufficiently large, the derivative 
\[
\frac{d}{ds} \Big|_{s=1} \area_{g}(\Sigma_{(s\xi,\lambda)})
\]
is negative when $|\xi| = 2 \sqrt 2$ and positive when $|\xi| = 5$. Using that the metric $g$ is rotationally symmetric, we see that the map 
\[
\xi \mapsto \area_{g} ( \Sigma_{(\xi,10^{j})})
\]
has a stable critical point (a local minimum) at some $\xi_{j} \in \R^3$ with $|\xi_j| \in (2 \sqrt 2, 5)$. In other words, $\Sigma_{(\xi_{j},10^{j})}$ is a ``far-off-center" stable constant mean sphere for $j$ sufficiently large. 


\begin{remark}
S. Brendle has already observed in Theorem 1.5 of \cite{Brendle:2013} that, as a consequence of the work by F. Pacard and X. Xu in \cite{PacardXu}, every rotationally symmetric Riemannian manifold whose scalar curvature has a strict local extremum contains small stable constant mean curvature spheres. 
\end{remark}


\section{Proof of Theorem \ref{theo:far-off-center-LS}} 

Consider 
\begin{align*} 
g_{ij}=(1+|x|^{-1})^{4}\delta_{ij} + \sigma_{ij}
\end{align*}
with 
\[
\partial_{I}\sigma_{ij} = O(|x|^{-2-|I|}) \qquad \text{ as } \qquad |x| \to \infty
\]
for all multi-indices $I$ of length $|I|\leq 7$. \\

Our proof is guided by the Lyapunov--Schmidt reduction and the related expansion for the reduced area functional as developed in \cite{Brendle-Eichmair:2014}. The goal is to extend these ideas to allow for $\xi\to\infty$. For a useful analysis in this regime, it is necessary to develop the expansion of the reduced area functional to a higher order than was necessary in \cite{Brendle-Eichmair:2014}, which turns out to be quite delicate. Our computations are also related and in part inspired by those for exact Schwarzschild in Appendix A of \cite{weingarten}. \\

We also note that part of our expansion for the reduced area functional $\area_{g}(\Sigma_{(\xi,\lambda)})$ follows, upon rescaling the chart at infinity by $\lambda  |\xi|$, from the work of S. Nardulli \cite{Nardulli:2009} or F.\ Pacard and X.\ Xu \cite{PacardXu}. 
The estimate for the error term in  \eqref{lsreduction} in e.g. \cite{PacardXu} is $O (\lambda^2 |\xi|^{-5})$ where we obtain $O (\lambda^{-1} |\xi|^{-6}) + O (|\xi|^{-7})$. Our stronger estimate is crucial for our applications here. \\

Let $\xi \in \R^3$ and $\lambda >0$ large. There is $r > 1$ with $r \sim \lambda$ and a smooth function $u_{(\xi, \lambda)}$ on the sphere $S_{r} ( \lambda \, \xi)$ that is perpendicular to constants and linear functions with respect to the Euclidean metric and such that the mean curvature with respect to $g$ of the Euclidean normal graph $\Sigma_{(\xi, \lambda)}$ of $u_{(\xi, \lambda)}$ --- as a function on $S_r( \lambda \, \xi)$ ---  is a linear combination of constants and linear functions and such that 
\[
\vol_{g}(\Sigma_{(\xi, \lambda)}) = 4\pi\lambda^{3}/3. 
\]
Moreover, 
\begin{align} \label{aux:estimateu0}
\sup_{S_r(\lambda \, \xi )} |u_{(\xi,\lambda)}| + \lambda \sup_{S_r(\lambda \, \xi )}  |\nabla u_{(\xi,\lambda)}| + \lambda^{2} \sup_{S_r(\lambda \, \xi )}  |\nabla^{2}u_{(\xi,\lambda)}| = O \left (1 / |\xi|\right).
\end{align}
This is a standard consequence of the implicit function theorem and elementary analysis; cf. Proposition 4 in \cite{Brendle-Eichmair:2014}. \\

We will improve estimate \eqref{aux:estimateu0} below. \\

It is convenient to abbreviate $a = \lambda \, \xi$. \\

We will frequently use the computations results listed in Appendix \ref{app:some-integral-expressions} in this section.


\subsection{Estimating \texorpdfstring{$\vol_{g}(B_{r}(a))$}{the volume in coordinate balls}}

Recall the following expansion for the determinant of a matrix
\begin{align*}
\sqrt{ \det(I+A)}  & = 1 + \frac 12 \tr A + \frac 1 8 (\tr A)^{2} - \frac 1 4 \tr A^{2} + O(|A|^{3}).
\end{align*}
Thus, we have  
\begin{align*}
& (1+|x|^{-1})^{6} \sqrt{\det(\delta_{ij} + (1+|x|^{-1})^{-4} \sigma_{ij})} \\
& \qquad = (1+|x|^{-1})^{6} \\
& \qquad \qquad+ \frac 12 (1+|x|^{-1})^{2} \tr  \sigma\\
& \qquad \qquad+  \frac 14 (1+|x|^{-1})^{-2} \left( \frac 1 2  (\tr  \sigma)^{2} - |\sigma|^{2} \right)\\
& \qquad \qquad+ O(|x|^{-6}).
\end{align*}
Repeating the computations in Proposition 17 of \cite{weingarten} (noting the dependence of the error on $r$), we find 
\begin{align*}
\int_{B_{r}(a)} (1+|x|^{-1})^{6} & = \frac{4\pi}{3} r^{3} (1+|a|^{-1})^{6} \Bigl( 1 + 3(1+|a|^{-1})^{-2}\frac{r^{2}}{|a|^{4}} + \frac{9}{7} \frac{r^{4}}{|a|^{6}}\Bigl) + O(r^{8}|a|^{-7})
\end{align*}
We now turn to the second term in the expansion of the volume form. We will write $\underline \sigma$ for $\sigma$ evaluated at $a$ (we will use the convention that if $\underline\sigma$ appears with a derivative, the derivative is taken \emph{and then} the quantity is evaluated at $a$). 

First, note that for $y \in B_{r}(0)$ with $x=a+y$,
\[
(1+|x|^{-1})^{2} = (1+|a|^{-1})^{2} + 2(1+|a|^{-1}) (|a+y|^{-1}-|a|^{-1}) + \underbrace{(|a+y|^{-1}-|a|^{-1})^{2}}_{=O(|a|^{-4}|y|^{2})}
\]
as well as
\[
|y+a|^{-1}-|a|^{-1} = - \frac{\bangle{a,y}}{|a|^{3}} - \frac 12 \frac{|a|^{2}|y|^{2}-3\bangle{a,y}^{2}}{|a|^{5}} + O(r^{3}|a|^{-4}).
\]
Finally, we have
\begin{align*}
\tr \sigma & = \tr  \underline \sigma + \nabla_{y}(\tr \underline \sigma) + \frac 12 \nabla^{2}_{y,y}(\tr \underline\sigma)\\
&\qquad + \frac 1 6 \nabla^{3}_{y,y,y}(\tr \underline\sigma) + \frac{1}{24} \nabla^{4}_{y,y,y,y}(\tr \underline\sigma) \\
& \qquad + O(|y|^{5}|x|^{-7}).
\end{align*}
We will frequently consider such Taylor expansions for expressions involving $\sigma$.

Combining the above expansions and using the expressions found in Appendix \ref{app:some-integral-expressions}, we have
\begin{align*}
\frac 12 \int_{B_{r}(a)} (1+|x|^{-1})^{2} \tr \sigma & = \frac 12 (1+|a|^{-1})^{2} \int_{B_{r}(a)} \tr \sigma\\
& \qquad+ (1+|a|^{-1})\int_{B_{r}(a)} (|a+y|^{-1}-|a|^{-1}) \tr  \sigma\\
& \qquad+ O(r^{5}|a|^{-6}) \\
& = \frac 12 (1+|a|^{-1})^{2} \int_{B_{r}} \tr \underline\sigma\\
&  \qquad+  \frac 14 (1+|a|^{-1})^{2} \int_{B_{r}} \nabla^{2}_{y,y}\tr \underline\sigma\\
&  \qquad+  \frac{1}{48} (1+|a|^{-1})^{2} \int_{B_{r}} \nabla^{4}_{y,y,y,y}\tr \underline\sigma\\
& \qquad- (1+|a|^{-1}) |a|^{-3}\int_{B_{r}} \bangle{a,y} \nabla_{y}\tr  \underline \sigma\\
& \qquad+ O(r^{5}|a|^{-6}) + O(r^{7}|a|^{-7}) \\
& = \frac{2\pi}{3} r^{3} (1+|a|^{-1})^{2}\tr \underline\sigma\\
& \qquad+ \frac{\pi}{15} r^{5} (1+|a|^{-1})^{2} \Delta(\tr \underline\sigma)\\
& \qquad+ \frac{\pi}{420} r^{7} \Delta( \Delta(\tr \underline\sigma))\\
& \qquad- \frac{4\pi}{15} (1+|a|^{-1})r^{5}|a|^{-3}\nabla_{a} (\tr \underline\sigma)\\
& \qquad+ O(r^{5}|a|^{-6}) + O(r^{7}|a|^{-7})
\end{align*}
Continuing on, we have that
\begin{align*}
\frac 14 \int_{B_{r}(a)} (1+|x|^{-1})^{-2} \Bigl( \frac 12 (\tr \sigma)^{2} - |\sigma|^{2}\Bigl) & = \frac{\pi}{3} r^{3} (1+|a|^{-1})^{-2} \Bigl( \frac 12 (\tr \underline \sigma)^{2} - |\underline \sigma|^{2}\Bigl)\\
& \qquad+ O(r^{5}|a|^{-6})
\end{align*}
Now, putting these terms together, we find that
\begin{align*}
\vol_{g}(B_{r}(a)) & = \frac{4\pi}{3} r^{3} (1+|a|^{-1})^{6} \Bigl( 1 + 3(1+|a|^{-1})^{-2}\frac{r^{2}}{|a|^{4}} + \frac{9}{7} \frac{r^{4}}{|a|^{6}}\Bigl)\\
& \qquad+ \frac{2\pi}{3} r^{3} (1+|a|^{-1})^{2}\tr \underline\sigma\\
& \qquad+ \frac{\pi}{15} r^{5} (1+|a|^{-1})^{2} \Delta(\tr \underline\sigma)\\
& \qquad+ \frac{\pi}{420} r^{7} \Delta( \Delta(\tr \underline\sigma))\\
& \qquad- \frac{4\pi}{15} (1+|a|^{-1}) r^{5}|a|^{-3}\nabla_{a} (\tr \underline\sigma)\\
& \qquad+ \frac{\pi}{3} r^{3} (1+|a|^{-1})^{-2} \Bigl( \frac 12 (\tr \underline \sigma)^{2} - |\underline \sigma|^{2}\Bigl)\\
& \qquad+ O(r^{5}|a|^{-6}) + O(r^{8}|a|^{-7}).
\end{align*}


\subsection{Estimating \texorpdfstring{$\area_{g}(S_{r}(a))$}{the area of coordinate spheres}} 
Using the above expansion, we have that the volume form of $S_{r}(a)$ becomes
\begin{align*}
d\mu_{g} & = (1+|x|^{-1})^{4}\sqrt{\det(\delta|_{S} + (1+|x|^{-1})^{-4}\sigma|_{S})}\\
& = (1+|x|^{-1})^{4} \\
&  \qquad+ \frac 12 \tr_{S}\sigma\\
&  \qquad+ \frac 1 4 (1+|x|^{-1})^{-4} \Bigl( \frac 12 (\tr_{S}\sigma)^{2} - |\sigma|_{S}|^{2}\Bigl)\\
&  \qquad+ O(|a|^{-6})\\
& = (1+|x|^{-1})^{4} \\
&  \qquad+ \frac 12 \tr \sigma -\frac 12 r^{-2}\sigma(y,y)\\
&  \qquad+ \frac 1 4 (1+|x|^{-1})^{-4} \Bigl( \frac 12 (\tr \sigma)^{2} - r^{-2}(\tr \sigma)\sigma(y,y)  - |\sigma|^{2} + 2 r^{-2}|\sigma(y,\, \cdot \, )|^{2} - \frac 12 r^{-4} \sigma(y,y)^{2}\Bigl)\\
&  \qquad+ O(|a|^{-6}).
\end{align*}
As in Proposition 17 of \cite{weingarten}, we have that
\[
\int_{S_{r}(a)} (1+|x|^{-1})^{4} = 4\pi r^{2}(1+|a|^{-1})^{4}\Bigl(1+2(1+|a|^{-1})^{-2} \frac{r^{2}}{|a|^{4}} + \frac 6 5 \frac{r^{4}}{|a|^{6}} \Bigl) + O(r^{7}|a|^{-7})
\]
We compute, using Appendix \ref{app:some-integral-expressions},
\begin{align*}
\frac 12 \int_{S_{r}(a)} \tr \sigma & = \frac 12 \int_{S_{r}} \tr \underline \sigma\\
&  \qquad+ \frac 14 \int_{S_{r}} \nabla^{2}_{y,y} \tr \underline\sigma\\
& \qquad + \frac{1}{48} \int_{S_{r}}\nabla^{4}_{y,y,y,y} \tr \underline\sigma\\
& \qquad + O(r^{6}|a|^{-8})\\
& = 2\pi r^{2}\tr \underline\sigma\\
& \qquad+ \frac{\pi}{3} r^{4}\Delta(\tr \underline\sigma)\\
& \qquad+ \frac{\pi}{60} r^{6}\Delta(\Delta(\tr \underline\sigma))\\
& \qquad+ O(r^{8}|a|^{-8})
\end{align*}
and
\begin{align*}
\frac 12 \int_{S_{r}(a)} r^{-2}\sigma(y,y) & = \frac 12 \int_{S_{r}} r^{-2}\underline\sigma(y,y)\\
&  \qquad+ \frac 14 \int_{S_{r}} r^{-2}\nabla^{2}_{y,y}\underline \sigma(y,y)\\
&  \qquad+ \frac{1}{48} \int_{S_{r}} r^{-2}\nabla^{4}_{y,y,y,y}\underline\sigma(y,y)\\
&  \qquad+ O(r^{8}|a|^{-8})\\
& = \frac{2\pi}{3} r^{2} \tr \underline \sigma\\
& \qquad+ \frac{\pi}{15} r^{4} \Delta(\tr \underline\sigma) + \frac{2\pi}{15} r^{4}\Div(\Div(\underline\sigma))\\
& \qquad+ \frac{\pi}{420} r^{6} \Delta(\Delta(\tr \underline\sigma)) + \frac{\pi}{105}\Delta(\Div(\Div(\underline\sigma))\\
& \qquad+ O(r^{8}|a|^{-8}).
\end{align*}
Putting these two expressions together, we find
\begin{align*}
\frac 12 \int_{S_{r}(a)} \tr \sigma - \frac 12 \int_{S_{r}(a)} r^{-2}\sigma(y,y) & =  2\pi r^{2}\tr \underline\sigma\\
& \qquad- \frac{2\pi}{3} r^{2} \tr \underline \sigma\\
& \qquad+ \frac{\pi}{3} r^{4}\Delta(\tr \underline\sigma)\\
& \qquad- \frac{\pi}{15} r^{4} \Delta(\tr \underline\sigma) - \frac{2\pi}{15} r^{4}\Div(\Div(\underline\sigma))\\
& \qquad+ \frac{\pi}{60} r^{6}\Delta(\Delta(\tr \underline\sigma))\\
& \qquad- \frac{\pi}{420} r^{6} \Delta(\Delta(\tr \underline\sigma)) - \frac{\pi}{105}\Delta(\Div(\Div(\underline\sigma))\\
& \qquad+ O(r^{8}|a|^{-8})\\
& =  \frac{4\pi}{3} r^{2}\tr \underline\sigma\\
& \qquad+ \frac{4\pi}{15} r^{4}\Delta(\tr \underline\sigma)  - \frac{2\pi}{15} r^{4}\Div(\Div(\underline\sigma))\\
& \qquad+ \frac{\pi}{70} r^{6}\Delta(\Delta(\tr \underline\sigma)) - \frac{\pi}{105}\Delta(\Div(\Div(\underline\sigma))\\
& \qquad+ O(r^{8}|a|^{-8}).
\end{align*}
Finally, we compute
\begin{align*}
& \frac 1 4 \int_{S_{r}(a)} (1+|x|^{-1})^{-4}\Bigl( \frac 12 (\tr\sigma)^{2} - r^{-2}(\tr\sigma)\sigma(y,y) - |\sigma|^{2} + 2r^{-2}|\sigma(y,\, \cdot \, )|^{2} - \frac 12 r^{-4}\sigma(y,y)^{2}\Bigl)\\
& = \frac 18 (1+|a|^{-1})^{-4}\int_{S_{r}} (\tr\underline\sigma)^{2}\\
&  \qquad- \frac 14 (1+|a|^{-1})^{-4} \int_{S_{r}} r^{-2}(\tr\underline\sigma)\underline\sigma(y,y)\\
&  \qquad- \frac 1 4 (1+|a|^{-1})^{-4} \int_{S_{r}} |\underline\sigma|^{2}\\
&  \qquad+ \frac 12 (1+|a|^{-1})^{-4} \int_{S_{r}} r^{-2} |\underline\sigma(y,\, \cdot \, )|^{2}\\
&  \qquad-\frac 18 (1+|a|^{-1})^{-4} \int_{S_{r}} r^{-4}\underline\sigma(y,y)^{2}\\
&  \qquad+ O(r^{4}|a|^{-6})\\
& = \frac{\pi}{2} r^{2} (1+|a|^{-1})^{-4} (\tr\underline\sigma)^{2}\\
& \qquad- \frac{\pi}{3} r^{2} (1+|a|^{-1})^{-4}(\tr\underline\sigma)^{2}\\
& \qquad- \pi r^{2} (1+|a|^{-1})^{-4} |\underline\sigma|^{2}\\
& \qquad+ \frac{2\pi}{3} r^{2} (1+|a|^{-1})^{-4} |\underline\sigma|^{2}\\
& \qquad- \frac{\pi}{30}r^{2} (1+|a|^{-1})^{-4}(\tr\underline\sigma)^{2} - \frac{\pi}{15} r^{2} (1+|a|^{-1})^{-4}|\underline\sigma|^{2}\\
& \qquad+ O(r^{4}|a|^{-6})\\
& = \frac{2\pi}{15} r^{2} (1+|a|^{-1})^{-4} (\tr\underline\sigma)^{2} - \frac{2\pi}{5} r^{2} (1+|a|^{-1})^{-4} |\underline \sigma|^{2} + O(r^{4}|a|^{-6})\\
& = - \frac{2\pi}{5} r^{2} (1+|a|^{-1})^{-4} |\underline {\mathring\sigma}|^{2} + O(r^{4}|a|^{-6}).
\end{align*}
Thus, putting this together, we find that
\begin{align*}
\area_{g}(S_{r}(a)) & = 4\pi r^{2}(1+|a|^{-1})^{4}\Bigl(1+2(1+|a|^{-1})^{-2} \frac{r^{2}}{|a|^{4}} + \frac 6 5 \frac{r^{4}}{|a|^{6}} \Bigl) \\
& \qquad+ \frac{4\pi}{3} r^{2}\tr \underline\sigma\\
& \qquad+ \frac{4\pi}{15} r^{4}\Delta(\tr \underline\sigma)  - \frac{2\pi}{15} r^{4}\Div(\Div(\underline\sigma))\\
& \qquad+ \frac{\pi}{70} r^{6}\Delta(\Delta(\tr \underline\sigma)) - \frac{\pi}{105}r^{6}\Delta(\Div(\Div(\underline\sigma))\\
& \qquad- \frac{2\pi}{5} r^{2} (1+|a|^{-1})^{-4} |\mathring {\underline \sigma}|^2 \\
& \qquad+ O(r^{4}|a|^{-6}) + O(r^{7}|a|^{-7})
\end{align*}


\subsection{Estimating \texorpdfstring{$\sF(S_{r}(a))$}{the brane functional of coordinate spheres}}
We define
\[
\sF(S_{r}(a)) = \area_{g}(S_{r}(a)) - 2r^{-1}(1+|a|^{-1})^{-2}\vol_{g}(S_{r}(a)).
\]
We then compute
\begin{align*}
\sF(S_{r}(a)) & =  4\pi r^{2}(1+|a|^{-1})^{4}\Bigl(1+2(1+|a|^{-1})^{-2} \frac{r^{2}}{|a|^{4}} + \frac 6 5 \frac{r^{4}}{|a|^{6}} \Bigl) \\
& \qquad+ \frac{4\pi}{3} r^{2}\tr \underline\sigma\\
& \qquad+ \frac{4\pi}{15} r^{4}\Delta(\tr \underline\sigma)  - \frac{2\pi}{15} r^{4}\Div(\Div(\underline\sigma))\\
& \qquad+ \frac{\pi}{70} r^{6}\Delta(\Delta(\tr \underline\sigma)) - \frac{\pi}{105}r^{6}\Delta(\Div(\Div(\underline\sigma))\\
& \qquad- \frac{2\pi}{5} r^{2} (1+|a|^{-1})^{-4} |\mathring {\underline \sigma}|^2 \\
& \qquad-\frac{8\pi}{3} r^{2} (1+|a|^{-1})^{6} \Bigl( 1 + 3(1+|a|^{-1})^{-2}\frac{r^{2}}{|a|^{4}} + \frac{9}{7} \frac{r^{4}}{|a|^{6}}\Bigl)\\
& \qquad- \frac{4\pi}{3} r^{2} \tr \underline\sigma\\
& \qquad- \frac{2\pi}{15} r^{4}  \Delta(\tr \underline\sigma)\\
& \qquad- \frac{\pi}{210} r^{5} \Delta( \Delta(\tr \underline\sigma))\\
& \qquad+ \frac{8\pi}{15} (1+|a|^{-1})^{-1} r^{5}|a|^{-3}\nabla_{a} (\tr \underline\sigma)\\
& \qquad- \frac{2\pi}{3} r^{3} (1+|a|^{-1})^{-4} \Bigl( \frac 12 (\tr \underline \sigma)^{2} - |\underline \sigma|^{2}\Bigl)\\
& \qquad+ O(r^{4}|a|^{-6}) + O(r^{7}|a|^{-7})\\
 & =  \frac{4\pi}{3} r^{2}(1+|a|^{-1})^{4} + \frac{48\pi}{35} \frac{r^{6}}{|a|^{6}}  \\
& \qquad+ \frac{2\pi}{15} r^{4}\Delta(\tr \underline\sigma)  - \frac{2\pi}{15} r^{4}\Div(\Div(\underline\sigma))\\
& \qquad+ \frac{\pi}{105} r^{6}\Delta(\Delta(\tr \underline\sigma)) - \frac{\pi}{105}r^{6}\Delta(\Div(\Div(\underline\sigma))\\
& \qquad- \frac{2\pi}{5} r^{2} (1+|a|^{-1})^{-4} |\mathring {\underline \sigma}|^2 \\
& \qquad- \frac{2\pi}{3} r^{2} (1+|a|^{-1})^{-4} \Bigl( \frac 12 (\tr \underline \sigma)^{2} - |\underline \sigma|^{2}\Bigl)\\
& \qquad+ \frac{8\pi}{15} r^{4}|a|^{-3}\nabla_{a} (\tr \underline\sigma)\\
& \qquad+ O(r^{4}|a|^{-6}) + O(r^{7}|a|^{-7}).
\end{align*}


\subsection{Estimating the mean curvature of \texorpdfstring{$S_{r}(a)$}{coordinate spheres}}

Consider 
\[
\hat g_{ij} = \bar g_{ij} + \hat \sigma_{ij} \qquad \text{ where } \qquad \hat \sigma_{ij} = (1+|x|^{-1})^{-4}\sigma_{ij}.
\]
By the computation in Lemma 7.4 of \cite{Huisken-Ilmanen:2001}, we have 
\begin{align*}
\hat H & =  H - r^{-1}\tr_{S} \hat \sigma + r^{-3}\hat \sigma(y,y) - r^{-1}\tr_{S}(\nabla_{\, \cdot \, }\hat\sigma)(y,\, \cdot \, ) + r^{-1} \frac 12 \tr_{S} \nabla_{y}\hat\sigma + O(r^{-1}|a|^{-4})\\
& = 2r^{-1} - r^{-1}\tr\hat\sigma + 2 r^{-3}\hat \sigma(y,y)- r^{-1} \Div(\hat \sigma)(y) + \frac 12 r^{-1} \nabla_{y}\tr \hat \sigma + \frac 12 r^{-3}\nabla_{y}\hat \sigma(y,y) + O(r^{-1}|a|^{-4})
\end{align*}
for the mean curvature of $S_r(a)$ with respect to $\hat g$. We recall the decomposition  
\[
a + y = x \in S_r(a)
\]
and that geometric quantities are computed with respect to the Euclidean background metric $\bar g$ unless noted otherwise. 
It follows that the mean curvature of $S_r(a)$ with respect to $g$ is given by 
\begin{align*}
H_g & = (1+|x|^{-1})^{-2} \hat H - 4(1+|x|^{-1})^{-3} |x|^{-3} \hat g(x,\hat \nu)\\
& =(1+|x|^{-1})^{-2} \hat H - 4r^{-1}(1+|x|^{-1})^{-3} |x|^{-3} \bangle{x,y} + O(|a|^{-4})\\
& = 2r^{-1}(1+|x|^{-1})^{-2} - 4r^{-1}(1+|x|^{-1})^{-3} |x|^{-3} \bangle{x,y} \\
& \qquad - r^{-1}(1+|x|^{-1})^{-6} \tr\sigma + 2 r^{-3}(1+|x|^{-1})^{-6} \sigma(y,y) \\
& \qquad- r^{-1}(1+|x|^{-1})^{-2} \Div(\hat \sigma)(y) + \frac 12 r^{-1}(1+|x|^{-1})^{-2}  \nabla_{y}\tr \hat \sigma + \frac 12 r^{-3}(1+|x|^{-1})^{-2}\nabla_{y}\hat \sigma(y,y)\\
& \qquad+ O(|a|^{-4})\\
& = 2r^{-1}(1+|x|^{-1})^{-2} - 4r^{-1}(1+|x|^{-1})^{-3} |x|^{-3} \bangle{x,y} \\
& \qquad - r^{-1}(1+|a|^{-1})^{-6} \tr\sigma + 2 r^{-3}(1+|a|^{-1})^{-6} \sigma(y,y) \\
& \qquad - r^{-1} \Div( \sigma)(y) + \frac 12 r^{-1} \nabla_{y}\tr  \sigma + \frac 12 r^{-3} \nabla_{y} \sigma(y,y)\\
& \qquad + O(r^{2}|a|^{-4}).
\end{align*}
Computing as in Lemma 18 of \cite{weingarten},
\begin{align*}
& 2r^{-1}(1+|x|^{-1})^{-2} - 4r^{-1}(1+|x|^{-1})^{-3} |x|^{-3} \bangle{x,y}\\
& \qquad= 2r^{-1}\bigl((1+|a|^{-1})^{-2} - ( |a|^{-3}|y|^{2}-3 |a|^{-5} \bangle{a,y}^{2} ) \bigl) \\
& \qquad \qquad +O(r^{2}|a|^{-4}).
\end{align*}
Thus,
\begin{align*}
H_g & = 2r^{-1}\bigl((1+|a|^{-1})^{-2} - (|a|^{-3}|y|^{2}-3 |a|^{-5}\bangle{a,y}^{2} ) \bigl)\\
& \qquad- r^{-1}(1+|a|^{-1})^{-6} \tr\sigma + 2 r^{-3}(1+|a|^{-1})^{-6} \sigma(y,y) \\
& \qquad- r^{-1} \Div( \sigma)(y) + \frac 12 r^{-1} \nabla_{y}\tr  \sigma + \frac 12 r^{-3} (\nabla_{y} \sigma)(y,y)\\
& \qquad+ O(r^{2}|a|^{-4}).
\end{align*}
Now, we consider the (Euclidean) projection of $H_g$ to $\Lambda_{2}$ and $\Lambda_{>2}$ where $\Lambda_{2}$ is the space of second eigenfunctions on $S_{r}$ and $\Lambda_{>2}$ is the $L^{2}(S_{r})$-orthogonal complement of $\Lambda_{0}\oplus \Lambda_{1}\oplus\Lambda_{2}$. 
\begin{align*}
\proj_{\Lambda_{2}} H_g & = - \frac{2}{r}\frac{|a|^{2}|y|^{2}-3\bangle{a,y}^{2}}{|a|^{5}} \\
& \qquad+ 2 r^{-3} (1+|a|^{-1})^{-6} \proj_{\Lambda_{2}}\underline\sigma(y,y)\\
& \qquad+ O(r^{2}|a|^{-4})\\
& = - - \frac{2}{r}\frac{|a|^{2}|y|^{2}-3\bangle{a,y}^{2}}{|a|^{5}} \\
& \qquad+ 2 r^{-3} (1+|a|^{-1})^{-6} \Bigl( \underline\sigma(y,y) - \frac 1 3 |y|^{2} \tr\underline\sigma \Bigl)\\
& \qquad+ O(r^{2}|a|^{-4}).
\end{align*}
For the higher eigenspaces, we will be content with the estimate
\[
\proj_{\Lambda_{>2}} H_g = O(|a|^{-3}) + O(r^{2}|a|^{-4}).
\]


\subsection{Estimates for \texorpdfstring{$u$}{}} Our goal here is to improve upon the initial estimate \eqref{aux:estimateu0}. \\

Let $t \in [0, 1]$. Consider the Euclidean graph over $S_r(a)$ of the function $t \, u$. The initial normal speed with respect to $g$ of this family can be computed as   
\[
w = u \, g(y/r,\nu_g). 
\]
Note that
\[
w = (1+O(|x|^{-1}))u
\]
up to and including second derivatives. We will give a more precise estimate later. Thus, the second variation of area implies that  
\[
\Delta^{S_{r}(a)}_g  w + (|h_g|_g^{2} + \Ric_g(\nu_g,\nu_g)) w = H_g- H_g^{\Sigma} + O(\lambda^{-3} |\xi|^{-2})
\]
where, as before, $H_g$ is the mean curvature of $S_r(a)$ with respect to $g$. It follows that 
\[
\Delta^{S_{r}(a)} u + 2r^{-2} u = H_g - H_g^{\Sigma} + O(\lambda^{-3} |\xi|^{-2}).
\]
Since 
\[
\proj_{\Lambda_{>1}} (H_g - H_g^{\Sigma}) = \proj_{\Lambda_{>1}} H_g = O(\lambda^{-3}|\xi|^{-2}) + O(\lambda^{-2}|\xi|^{-3}),
\] 
we obtain that  
\[
\sup_{S_r(\lambda \, \xi )} |u_{(\xi,\lambda)}| + \lambda \sup_{S_r(\lambda \, \xi )}  |\nabla u_{(\xi,\lambda)}| + \lambda^{2} \sup_{S_r(\lambda \, \xi )}  |\nabla^{2}u_{(\xi,\lambda)}|  =  O(\lambda^{-1}|\xi|^{-2}) + O(|\xi|^{-3}).
\]
This allows us to improve the coarse estimate above to 
\[
\Delta_g^{S_{r}(a)} w + (|h_g|_g^{2} + \Ric_g(\nu_g,\nu_g)) w = H_g - H_g^{\Sigma} + O(\lambda^{-5}|\xi|^{-4}) + O(\lambda^{-3} |\xi|^{-6}) .
\]
At this point, we can improve our earlier estimate for $w$ to
\[
w = \left((1+|x|^{-1})^{2} + O(|x|^{-2})\right) u
\]
up to and including second derivatives. Thus
\begin{align*}
\Delta_g^{S_{r}(a)} w & = (1+|a|^{-1})^{-2}\Delta^{S_{r}(a)} u + O(\lambda^{-4}|\xi|^{-4}) + O(\lambda^{-3}|\xi|^{-5}).
\end{align*}
Continuing on, we have that 
\[
|h_g|_g^{2} = 2r^{-2}(1+|a|^{-1})^{-4} + O(\lambda^{-4}|\xi|^{-2}) 
\]
and
\[
\Ric_g(\nu_g,\nu_g) = O(\lambda^{-3}|\xi|^{-3}).
\]
Putting these estimates together, we find that
\[
(1+|a|^{-1})^{-2}\Delta^{S_{r}(a)} u + 2r^{-2}(1+|a|^{-1})^{-2} u = H_g-H_g^{\Sigma} + O(\lambda^{-4}|\xi|^{-4}) + O(\lambda^{-3}|\xi|^{-5}).
\]
Hence,
\begin{align*}
& \Delta^{S_{r}(a)}\proj_{\Lambda_{2}} u + 2r^{-2}\proj_{\Lambda_{2}} u \\
& \qquad =  \proj_{\Lambda_{2}} (   \Delta^{S_{r}(a)}u + 2r^{-2} u) \\
& \qquad= (1 + |a|^{-1})^{2}  \proj_{\Lambda_{2}}H_g +  O(\lambda^{-4}|\xi|^{-4}) + O(\lambda^{-3}|\xi|^{-5})\\
& \qquad=  \frac{2}{r^3} \frac{1}{(1 + |a|^{-1})^{4}}\Bigl(\underline\sigma(y,y) -\frac 1 3 |y|^{2}\tr\underline\sigma\Bigl) -\frac{2}{r}\frac{|a|^{2}|y|^{2}-3\bangle{a,y}^{2}}{|a|^{5}}   + O(\lambda^{-2}|\xi|^{-4}).
\end{align*}
This implies that
\[
\proj_{\Lambda_{2}} u = -\frac{1}{2 \, r} \frac{1}{(1+|a|^{-1})^{4}}\Bigl(\underline\sigma(y,y) -\frac 1 3 |y|^{2}\tr\underline\sigma\Bigl) + \frac r 2 \frac{|a|^{2}|y|^{2}-3\bangle{a,y}^{2}}{|a|^{5}}   + O(|\xi|^{-4})
\]
together with two derivatives. Note that in particular 
\[
\proj_{\Lambda_{2}} u = O(\lambda^{-1}|\xi|^{-2}) + O(|\xi|^{-3})
\]
along with two derivatives. The above expression also implies that
\[
\proj_{\Lambda_{>2}} u = O(\lambda^{-1} |\xi|^{-3}) + O(|\xi|^{-4})
\]
with two derivatives. 


\subsection{Estimating \texorpdfstring{$\sF(\Sigma)$}{the brane functional of the perturbed surface}}

We have that
\begin{align*}
\sF(\Sigma) & = \sF(S_{r}(a)) + \int_{S_{r}(a)} (H_g - 2r^{-1}(1+|a|^{-1})^{-2}) w \, d\mu_{g}\\
& \qquad + \frac 12 \int_{S_{r}(a)} H_g (H_g - 2r^{-1}(1+|a|^{-1})^{-2}) w^{2} d\mu_{g}\\
& \qquad - \frac 12 \int_{S_{r}(a)} (\Delta_g^{S_{r}(a)} w + (|h_g|_g^{2} + \Ric_g(\nu_g,\nu_g)) w)w\, d\mu_{g}\\
& \qquad+ O(\lambda^{-4}|\xi|^{-6}) + O(\lambda^{-1}|\xi|^{-9})
\end{align*}
Recall that 
\[
w = g(y/r,\nu_g) = (1+|x|^{-1})^{2}\Bigl(1 + \frac 12 (1+|x|^{-1})^{-4}r^{-2}\sigma(y,y) \Bigl)+ O(|x|^{-4}).
\]
We have seen above that
\[
d\mu_{g} = (1+|x|^{-1})^{4} \Bigl(1 + \frac 12 (1+|x|^{-1})^{-4}\tr\sigma - \frac 12  (1+|x|^{-1})^{-4} r^{-2}\sigma(y,y)  \Bigl) d\mu_{\bar g} + O(|x|^{-4}).
\]
We begin with the first term.  
\begin{align*}
\int_{S_{r}(a)} (H_g-2r^{-1}(1+|a|^{-1})^{-2}) \, w \, d\mu_{g} & = \int_{S_{r}(a)} (H_g-2r^{-1}(1+|a|^{-1})^{-2}) \, u \, (1+|x|^{-1})^{6}\\
& \qquad+ O(\lambda^{-4}|\xi|^{-6}) + O(\lambda^{-3}|\xi|^{-7})\\
& = (1+|a|^{-1})^{6}\int_{S_{r}(a)} (H_g-2r^{-1}(1+|a|^{-1})^{-2}) \, u \\
& \qquad+ O(\lambda^{-3}|\xi|^{-6}) + O(\lambda^{-2}|\xi|^{-7})\\
& = - 2 r^{-1}(1+|a|^{-1})^{6} \int_{S_{r}(a)} \Bigl(\frac{|a|^{2}|y|^{2}-3\bangle{a,y}^{2}}{|a|^{5}}\Bigl) \, u\\
& \qquad+2 r^{-3} \int_{S_{r}(a)} u\, \underline\sigma(y,y)\\
& \qquad+ O(\lambda^{-2}|\xi|^{-6}) + O(\lambda^{-1}|\xi|^{-7})\\
& = - \int_{S_{r}} \Bigl(\frac{|a|^{2}|y|^{2}-3\bangle{a,y}^{2}}{|a|^{5}}\Bigl)^{2} \\
& \qquad+ r^{-2} (1+|a|^{-1})^{2} \int_{S_{r}} \Bigl(\frac{|a|^{2}|y|^{2}-3\bangle{a,y}^{2}}{|a|^{5}}\Bigl) \, \underline\sigma(y,y)\\
& \qquad- r^{-4} (1+|a|^{-1})^{-4} \int_{S_{r}}\Bigl(\underline\sigma(y,y)-\frac 13 |y|^{2}\tr\underline\sigma\Bigl)^{2}\\
& \qquad+ r^{-2} \int_{S_{r}} \Bigl(\frac{|a|^{2}|y|^{2}-3\bangle{a,y}^{2}}{|a|^{5}}\Bigl) \, \underline\sigma(y,y)\\
& \qquad+ O(\lambda^{-1}|\xi|^{-6}) + O(|\xi|^{-7})\\
& = - \int_{S_{r}} \Bigl(\frac{|a|^{2}|y|^{2}-3\bangle{a,y}^{2}}{|a|^{5}}\Bigl)^{2} \\
& \qquad+ \frac{2}{r^2}\int_{S_{r}} \Bigl(\frac{|a|^{2}|y|^{2}-3\bangle{a,y}^{2}}{|a|^{5}}\Bigl) \, \underline\sigma(y,y)\\
& \qquad \frac{1}{r^4} \frac{1}{(1+|a|^{-1})^{4}} \int_{S_{r}}\Bigl(\underline\sigma(y,y)-\frac 13 |y|^{2}\tr\underline\sigma\Bigl)^{2}\\
& \qquad+ O(\lambda^{-1}|\xi|^{-6}) + O(|\xi|^{-7}).
\end{align*}
The second term satisfies 
\begin{align*}
\frac 12 \int_{S_{r}(a)} H_g (H_g - 2r^{-1}(1+|a|^{-1})^{-2}) w^{2} d\mu_{g} & = O(\lambda^{-4}|\xi|^{-6}) + O(\lambda^{-1}|\xi|^{-9}).
\end{align*}
Finally, the last term satisfies 
\begin{align*}
& - \frac 12 \int_{S_{r}(a)} (\Delta_g^{S_{r}(a)} w + (|h_g|_g^{2} + \Ric_g(\nu_g,\nu_g)) w)w\, d\mu_{g}\\
& \qquad=  - \frac 12 (1+|a|^{-1})^{4}  \int_{S_{r}(a)} (\Delta^{S_{r}} u + 2 \, r^{-2}u)u\\
& \qquad\qquad + O(\lambda^{-3}|\xi|^{-6}) + O(\lambda^{-1}|\xi|^{-8})\\
& \qquad=  2r^{-2} (1+|a|^{-1})^{4}  \int_{S_{r}(a)} (\proj_{\Lambda_{2}}u)^{2}\\
& \qquad\qquad + O(\lambda^{-2}|\xi|^{-6}) + O(|\xi|^{-8})\\
& \qquad=   \frac{1}{2 \, r^4} \frac{1}{ (1+|a|^{-1})^{4}}  \int_{S_{r}} \Bigl(\underline\sigma(y,y) - \frac 1 3 |y|^{2} \tr\underline\sigma \Bigl)^{2}\\
& \qquad\qquad -  \frac{1}{r^2} \int_{S_{r}} \Bigl(\frac{|a|^{2}|y|^{2}-3\bangle{a,y}^{2}}{|a|^{5}}\Bigl)\underline\sigma(y,y) \\
&\qquad \qquad +\frac 12 \int_{S_{r}} \Bigl(\frac{|a|^{2}|y|^{2}-3\bangle{a,y}^{2}}{|a|^{5}}\Bigl)^{2} \\
&  \qquad\qquad+ O(\lambda^{-1}|\xi|^{-6}) + O(|\xi|^{-7}).
\end{align*}
Putting this together, we find that 
\begin{align*}
\sF(\Sigma) & = \sF(S_{r}(a)) \\
& \qquad- \int_{S_{r}} \Bigl(\frac{|a|^{2}|y|^{2}-3\bangle{a,y}^{2}}{|a|^{5}}\Bigl)^{2} \\
& \qquad+ 2r^{-2}\int_{S_{r}} \Bigl(\frac{|a|^{2}|y|^{2}-3\bangle{a,y}^{2}}{|a|^{5}}\Bigl)\underline\sigma(y,y)\\
& \qquad- r^{-4} (1+|a|^{-1})^{-4} \int_{S_{r}}\Bigl(\underline\sigma(y,y)-\frac 13 |y|^{2}\tr\underline\sigma\Bigl)^{2}\\
& \qquad+ \frac 12 r^{-4} (1+|a|^{-1})^{-4}  \int_{S_{r}} \Bigl(\underline\sigma(y,y) - \frac 1 3 |y|^{2} \tr\underline\sigma \Bigl)^{2}\\
& \qquad-  r^{-2}  \int_{S_{r}} \Bigl(\frac{|a|^{2}|y|^{2}-3\bangle{a,y}^{2}}{|a|^{5}}\Bigl)\underline\sigma(y,y) \\
& \qquad+ \frac 12 \int_{S_{r}} \Bigl(\frac{|a|^{2}|y|^{2}-3\bangle{a,y}^{2}}{|a|^{5}}\Bigl)^{2} \\
& \qquad+ O(\lambda^{-1}|\xi|^{-6}) + O(|\xi|^{-7})\\
& = \sF(S_{r}(a)) \\
& \qquad- \frac 12 \int_{S_{r}} \Bigl(\frac{|a|^{2}|y|^{2}-3\bangle{a,y}^{2}}{|a|^{5}}\Bigl)^{2} \\
& \qquad+ r^{-2}\int_{S_{r}} \Bigl(\frac{|a|^{2}|y|^{2}-3\bangle{a,y}^{2}}{|a|^{5}}\Bigl)\underline\sigma(y,y)\\
& \qquad- \frac 12 r^{-4} (1+|a|^{-1})^{-4} \int_{S_{r}}\Bigl(\underline\sigma(y,y)-\frac 13 |y|^{2}\tr\underline\sigma\Bigl)^{2}\\
& \qquad+ O(\lambda^{-1}|\xi|^{-6}) + O(|\xi|^{-7})\\
\intertext{We now use the expansions given in Appendix \ref{sec:someintegrals}.}
& =  \sF(S_{r}(a)) \\
& \qquad- \frac{8\pi}{5} |\xi|^{-6}\\
& \qquad+ \frac{8\pi}{15} \frac{r^{4}}{|a|^{3}}\left(\tr\underline\sigma - 3|a|^{-1}\underline\sigma(a,a)\right)\\
& \qquad- \frac{4\pi}{15} r^{2}(1+|a|^{-1})^{-4} |\mathring{\underline\sigma}|^{2}\\
& \qquad+ O(\lambda^{-1}|\xi|^{-6}) + O(|\xi|^{-7})\\
& =  \frac{4\pi}{3} r^{2}(1+|a|^{-1})^{4} + \frac{48\pi}{35} \frac{r^{6}}{|a|^{6}}  \\
& \qquad+ \frac{2\pi}{15} r^{4}\Delta(\tr \underline\sigma)  - \frac{2\pi}{15} r^{4}\Div(\Div(\underline\sigma))\\
& \qquad+ \frac{\pi}{105} r^{6}\Delta(\Delta(\tr \underline\sigma)) - \frac{\pi}{105}r^{6}\Delta(\Div(\Div(\underline\sigma))\\
&\qquad - \frac{2\pi}{5} r^{2} (1+|a|^{-1})^{-4} |\mathring {\underline \sigma}|^2 \\
& \qquad- \frac{2\pi}{3} r^{2} (1+|a|^{-1})^{-4} \Bigl( \frac 12 (\tr \underline \sigma)^{2} - |\underline \sigma|^{2}\Bigl)\\
& \qquad+ \frac{8\pi}{15} r^{4}|a|^{-3}\nabla_{a} (\tr \underline\sigma)\\
&\qquad - \frac{8\pi}{5} |\xi|^{-6}\\
& \qquad+ \frac{8\pi}{15} \frac{r^{4}}{|a|^{3}}\left(\tr\underline\sigma - 3|a|^{-1}\underline\sigma(a,a)\right)\\
&\qquad - \frac{4\pi}{15} r^{2}(1+|a|^{-1})^{-4} |\mathring{\underline\sigma}|^{2}\\
&\qquad + O(\lambda^{-1}|\xi|^{-6}) + O(|\xi|^{-7})\\
& =  \frac{4\pi}{3} r^{2}(1+|a|^{-1})^{4} - \frac{8\pi}{35} |\xi|^{-6}  \\
& \qquad+ \frac{2\pi}{15} r^{4}\Delta(\tr \underline\sigma)  - \frac{2\pi}{15} r^{4}\Div(\Div(\underline\sigma))\\
& \qquad+ \frac{\pi}{105} r^{6}\Delta(\Delta(\tr \underline\sigma)) - \frac{\pi}{105}r^{6}\Delta(\Div(\Div(\underline\sigma))\\
& \qquad- \frac{2\pi}{3} r^{2} (1+|a|^{-1})^{-4} |\mathring {\underline \sigma}|^2 \\
&\qquad - \frac{2\pi}{3} r^{2} (1+|a|^{-1})^{-4} \Bigl( \frac 12 (\tr \underline \sigma)^{2} - |\underline \sigma|^{2}\Bigl)\\
& \qquad+ \frac{8\pi}{15} r^{4}|a|^{-3}\nabla_{a} (\tr \underline\sigma)\\
& \qquad+ \frac{8\pi}{15} \frac{r^{4}}{|a|^{3}}\left(\tr\underline\sigma - 3|a|^{-1}\underline\sigma(a,a)\right)\\
& \qquad+ O(\lambda^{-1}|\xi|^{-6}) + O(|\xi|^{-7}).
\end{align*}
Using that $\vol_{g}(\Omega) = \frac{4\pi}{3}\lambda^{3}$, we obtain 
\begin{align*}
\area_{g}(\Sigma) & =  \frac{4\pi}{3} r^{2}(1+|a|^{-1})^{4} + \frac{8\pi}{3} \lambda^{3} r^{-1} (1+|a|^{-1})^{-2} - \frac{8\pi}{35}|\xi|^{-6}  \\
& \qquad+ \frac{2\pi}{15} r^{4}\Delta(\tr \underline\sigma)  - \frac{2\pi}{15} r^{4}\Div(\Div(\underline\sigma))\\
& \qquad+ \frac{\pi}{105} r^{6}\Delta(\Delta(\tr \underline\sigma)) - \frac{\pi}{105}r^{6}\Delta(\Div(\Div(\underline\sigma))\\
& \qquad- \frac{2\pi}{3} r^{2} (1+|a|^{-1})^{-4} |\mathring {\underline \sigma}|^2 \\
& \qquad- \frac{2\pi}{3} r^{2} (1+|a|^{-1})^{-4} \Bigl( \frac 12 (\tr \underline \sigma)^{2} - |\underline \sigma|^{2}\Bigl)\\
& \qquad+ \frac{8\pi}{15} r^{4}|a|^{-3}\nabla_{a} (\tr \underline\sigma)\\
& \qquad+ \frac{8\pi}{15} \frac{r^{4}}{|a|^{3}}\left(\tr\underline\sigma - 3|a|^{-2}\underline\sigma(a,a)\right)\\
& \qquad+ O(\lambda^{-1}|\xi|^{-6}) + O(|\xi|^{-7}).
\end{align*}


\subsection{Estimating \texorpdfstring{$r$}{the radius}} 
We now use the expansion
\begin{align*}
\vol_{g}(\Omega) & = \vol_{g}(B_{r}(a)) + \int_{S_{r}(a)} w\, d\mu_{g} + \frac 12 \int_{S_{r}(a)} H_g w^{2}d\mu_{g} \\
& \qquad+ O( \lambda^{-3}|\xi|^{-6}) + O(|\xi|^{-9})
\end{align*}
to relate $\lambda$ and $r$. Note that because $u$ is orthogonal to constants and to linear functions,
\[
\int_{S_{r}(a)} w \, d\mu_{g} = O(|\xi|^{-5}) + O(\lambda^{-1}|\xi|^{-6})
\]
and 
\[
\frac 12 \int_{S_{r}(a)} H_g w^{2}d\mu_{g} = O(\lambda^{-1}|\xi|^{-4}) + O(\lambda|\xi|^{-6}).
\]
Hence, using the expression for $\vol_{g}(B_{r}(a))$ obtained previously, we find that 
\begin{align*}
\frac{4\pi}{3} \lambda^{3} & = \frac{4\pi}{3} r^{3}(1+|a|^{-1})^{6} + \frac{2\pi}{3} r^{3} (1+|a|^{-1})^{2} \tr\underline\sigma + O(\lambda |\xi|^{-4})\\
& = \frac{4\pi}{3} r^{3}(1+|a|^{-1})^{6} \Big( 1 + \frac 12 (1+|a|^{-1})^{-4} \tr\underline\sigma + O(\lambda^{-2}|\xi|^{-4}) \Big) .
\end{align*}

It is convenient to write
\[
\lambda^{3} = r^{3}(1+|a|^{-1})^{6} (1+\psi)
\]
for 
\[
\psi = \frac 12 (1+|a|^{-1})^{-4} \tr\underline\sigma + O(\lambda^{-2}|\xi|^{-4}) = O(\lambda^{-2}|\xi|^{-2})
\]
We now estimate the first line in the expansion for $\area_{g}(\Sigma)$ obtained above. 
\begin{align*}
& \frac{4\pi}{3} r^{2}(1+|a|^{-1})^{4} + \frac{8\pi}{3} \lambda^{3}r^{-1}(1+|a|^{-1})^{-2} \\
& \qquad = \frac{4\pi}{3} r^{2}(1+|a|^{-1})^{4} + \frac{8\pi}{3}r^{2} (1+|a|^{-1})^{4}(1+\psi)\\
&\qquad =  4\pi r^{2}(1+|a|^{-1})^{4} \left( 1 + \frac 2 3 \psi \right) \\
& \qquad=  4\pi r^{2}(1+|a|^{-1})^{4} \left( 1 +  \psi \right)^{\frac 2 3} +  \frac{4\pi}{9}r^{2} (1+|a|^{-1})^{4} \psi^{2} + O(r^{2}\psi^{3}) \\
& \qquad= 4\pi \lambda^{2} + \frac{\pi}{9} r^{2}(1+|a|^{-1})^{-4} (\tr\underline\sigma)^{2} + O(\lambda^{-2}|\xi|^{-6}).
\end{align*}


\subsection{Concluding the estimate for \texorpdfstring{$\area_{g}(\Sigma)$}{area}}
Combining the previous two subsections, we conclude that
\begin{align*}
\area_{g}(\Sigma) & = 4\pi \lambda^{2} - \frac{8\pi}{35}|\xi|^{-6}  \\
& \qquad + \frac{\pi}{9} r^{2} (1+|a|^{-1})^{-4} (\tr\underline\sigma)^{2} \\
& \qquad+ \frac{2\pi}{15} r^{4}\Delta(\tr \underline\sigma)  - \frac{2\pi}{15} r^{4}\Div(\Div(\underline\sigma))\\
& \qquad+ \frac{\pi}{105} r^{6}\Delta(\Delta(\tr \underline\sigma)) - \frac{\pi}{105}r^{6}\Delta(\Div(\Div(\underline\sigma))\\
& \qquad- \frac{2\pi}{3} r^{2} (1+|a|^{-1})^{-4} |\mathring {\underline \sigma}|^2 \\
& \qquad- \frac{2\pi}{3} r^{2} (1+|a|^{-1})^{-4} \Bigl( \frac 12 (\tr \underline \sigma)^{2} - |\underline \sigma|^{2}\Bigl)\\
& \qquad+ \frac{8\pi}{15} r^{4}|a|^{-3}\nabla_{a} (\tr \underline\sigma)\\
& \qquad+ \frac{8\pi}{15} \frac{r^{4}}{|a|^{3}}\left(\tr\underline\sigma - 3|a|^{-2}\underline\sigma(a,a)\right)\\
& \qquad+ O(\lambda^{-1}|\xi|^{-6}) + O(|\xi|^{-7})\\
& = 4\pi \lambda^{2} - \frac{8\pi}{35}|\xi|^{-6}  \\
& \qquad+ \frac{2\pi}{15} r^{4}\Delta(\tr \underline\sigma)  - \frac{2\pi}{15} r^{4}\Div(\Div(\underline\sigma))\\
& \qquad+ \frac{\pi}{105} r^{6}\Delta(\Delta(\tr \underline\sigma)) - \frac{\pi}{105}r^{6}\Delta(\Div(\Div(\underline\sigma))\\
& \qquad+ \frac{8\pi}{15} r^{4}|a|^{-3}\nabla_{a} (\tr \underline\sigma)\\
& \qquad+ \frac{8\pi}{15} \frac{r^{4}}{|a|^{3}}\left(\tr\underline\sigma - 3|a|^{-2}\underline\sigma(a,a)\right)\\
& \qquad+ O(\lambda^{-1}|\xi|^{-6}) + O(|\xi|^{-7})\\
& = 4\pi \lambda^{2} - \frac{8\pi}{35}|\xi|^{-6}  \\
&\qquad + \frac{2\pi}{15} \lambda^{4}(1+|a|^{-1})^{-8} \left(\Delta(\tr \underline\sigma)  - \Div(\Div(\underline\sigma))\right)\\
& \qquad+ \frac{\pi}{105} \lambda^{6} \left( \Delta(\Delta(\tr \underline\sigma)) -  \Delta(\Div(\Div(\underline\sigma))\right)\\
& \qquad+ \frac{8\pi}{15} \lambda^{4}|a|^{-3}\nabla_{a} (\tr \underline\sigma)\\
& \qquad+ \frac{8\pi}{15} \frac{\lambda^{4}}{|a|^{3}}\left(\tr\underline\sigma - 3|a|^{-2}\underline\sigma(a,a)\right)\\
& \qquad+ O(\lambda^{-1}|\xi|^{-6}) + O(|\xi|^{-7}).
\end{align*} 


\subsection{Estimating \texorpdfstring{$R$ and $\Delta_g R$}{scalar curvature}} 
We now relate the previous expression to the scalar curvature $R$ of $(M, g)$. As with mean curvature, we first consider 
\[
\hat g_{ij} = \bar g_{ij} + \hat \sigma_{ij} \qquad \text{ where } \qquad \hat \sigma_{ij} = (1+|x|^{-1})^{-4}\sigma_{ij}.
\]
Then,
\[
R_{\hat g} = \Div\Div\hat \sigma  - \Delta \tr\hat \sigma + O(|x|^{-6}).
\]
Note that
\[
\Div \hat \sigma = (1+|x|^{-1})^{-4} \Div\sigma + 4 (1+|x|^{-1})^{-5} |x|^{-3}\sigma(x,\, \cdot \, ).
\]
Thus, we find that
\begin{align*}
\Div\Div\hat\sigma & = (1+|x|^{-1})^{-4}\Div\Div\sigma + 4(1+|x|^{-1})^{-5}|x|^{-3}\Div \sigma(x)\\
& \qquad+ 20(1+|x|^{-1})^{-6} |x|^{-6} \sigma(x,x) - 12 (1+|x|^{-1})^{-5} |x|^{-5} \sigma(x,x) \\
& \qquad + 4(1+|x|^{-1})^{-5}|x|^{-3}\Div\sigma(x) + 4(1+|x|^{-1})^{-5}|x|^{-3} \tr\sigma\\
& = (1+|x|^{-1})^{-4}\Div\Div\sigma + 8 \Div \sigma(x)\\
& \qquad+ 4(1+|x|^{-1})^{-5}|x|^{-3} ( \tr\sigma - 3 |x|^{-2}\sigma(x,x)) \\
& \qquad+ O(|x|^{-6}).
\end{align*}
Similarly,
\begin{align*}
\Delta \tr\hat \sigma & = \Delta \left( (1+|x|^{-1})^{-4} \tr\sigma\right)\\
& = (1+|x|^{-1})^{-4}\Delta\tr\sigma \\
& \qquad+ 8(1+|x|^{-1})^{-5} |x|^{-3} \nabla_{x}\tr\sigma\\
& \qquad+ (\tr\sigma )\, \Delta (1+|x|^{-1})^{-4}\\
& = (1+|x|^{-1})^{-4} \, \Delta\tr\sigma \\
& \qquad+ 8 |x|^{-3} \nabla_{x}\tr\sigma\\
& \qquad+ O(|x|^{-6}).
\end{align*}
Thus, we find that 
\begin{align*}
R_{\hat g} & = (1+|x|^{-1})^{-4} \left( \Div\Div\sigma - \Delta \tr\sigma \right) \\
&  \qquad + 4 |x|^{-3} (\tr \sigma - 3 |x|^{-2} \sigma(x,x))\\
& \qquad+ 8 |x|^{-3}\Div(\sigma)(x) \\
& \qquad - 8 |x|^{-3} \nabla_{x}\tr\sigma\\
& \qquad + O(|x|^{-6}).
\end{align*}
It follows that 
\begin{align*}
R & =  - 8 (1 + |x|^{-1})^{-5} \Delta_{\hat g} |x|^{-1} + (1 + |x|^{-1})^{-4} R_{\hat g} \\
& = - 8 (1 + |x|^{-1})^{-5} \Delta_{\hat g} |x|^{-1} \\
& \qquad + (1 + |x|^{-1})^{-8}\left( \Div(\Div(\sigma)) - \Delta \tr\sigma\right) \\
&  \qquad + 4 |x|^{-3} ( \tr \sigma - 3 |x|^{-2} \sigma(x,x))\\
& \qquad+ 8 |x|^{-3}\Div(\sigma)(x) \\
& \qquad - 8 |x|^{-3} \nabla_{x}\tr\sigma\\
& \qquad + O(|x|^{-6}).
\end{align*}
Thus, it remains to estimate $\Delta_{\hat g} |x|^{-1}$. We have that 
\[
\sqrt{\det \hat g} = \sqrt{\det(\delta_{ij} + \hat \sigma_{ij})} = 1 + \frac 12 \tr \hat \sigma + O(|x|^{-4})
\] 
and
\[
\hat g^{ij} = \delta^{ij} - \hat \sigma^{ij} + O(|x|^{-4}).
\]
Thus,
\begin{align*}
\Delta_{\hat g} |x|^{-1} & = - 3 |x|^{-5} \sigma(x,x) \\
& \qquad+  |x|^{-3} \tr\sigma\\
& \qquad+   |x|^{-3}  \Div \sigma(x)\\
&\qquad - \frac 12 |x|^{-3}\nabla_{x}\tr\sigma \\
& \qquad+ O(|x|^{-6}).
\end{align*}
Thus, we find that 
\begin{align*}
R & = (1+|x|^{-1})^{-8} \left( \Div(\Div(\sigma)) - \Delta \tr\sigma\right) \\
& \qquad- 4|x|^{-3} \left( \tr\sigma - 3 |x|^{-2} \sigma (x, x) \right) \\
& \qquad-4   |x|^{-3}  \nabla_x \tr \sigma\\
& \qquad+ O(|x|^{-6}).
\end{align*}
Similarly, 
\begin{align*}
\Delta R & = \Delta  \left( \Div(\Div(\sigma)) - \Delta \tr\sigma \right) \\
& \qquad  + O (|x|^{-7}).
\end{align*}


\subsection{Reduced area-functional}

We finally obtain that, for $\xi \in \R^3$ and $\lambda > 0$ large,  
\begin{align*} 
(\xi, \lambda) \mapsto \area_{g}(\Sigma_{(\xi, \lambda)}) 
 = 4\pi \lambda^{2} - \frac{2\pi}{15} \lambda^{4}\underline R  - \frac{\pi}{105} \lambda^{6} \Delta\underline R  - \frac{8\pi}{35} |\xi|^{-6}  + O(\lambda^{-1} |\xi|^{-6})  
 + O(|\xi|^{-7})
\end{align*} 
where $R$ is the scalar curvature of $(M, g)$, and where 
\[
\underline R = R (\lambda \, \xi) \qquad \text{ and } \qquad \Delta \underline R = (\Delta R ) (\lambda \, \xi).
\]
The Laplacian is computed with respect to the Euclidean background metric. This is \eqref{lsreduction}.

We also record here the first radial derivative 
\begin{align} \label{lsradial}
\frac{d}{ds} \Big|_{s= 1} \area_g (\Sigma_{(s \xi, \lambda}) & = - \frac{2 \pi}{15} \lambda^5 |\xi| \partial_r \underline R - \frac{\pi}{105} \lambda^7  |\xi| \partial_r \Delta \underline R + \frac{48 \pi}{35} |\xi|^{-6} \nonumber \\
& \qquad  + O (\lambda^{-1} |\xi|^{-6}) + O (|\xi|^{-7}) \nonumber \\
& =  \frac{\pi}{105} \left( - 14 \lambda^5 |\xi| \partial_r \underline R - \lambda^7 |\xi| \partial_r \Delta \underline R + 144 \,  |\xi|^{-6}  \right) \\
&\qquad + O (\lambda^{-1} |\xi|^{-6}) + O (|\xi|^{-7}). \nonumber 
\end{align}
where the underscore indicates evaluation at $\lambda \, \xi$ after all derivatives are taken.

This completes the proof of Theorem \ref{theo:far-off-center-LS}.


\section{Proof of Corollary \ref{cor:veryfar-concave}}

We assume that $(M,g)$ is $C^{6}$-asymptotically Schwarzschild in the sense that 
\[
g_{ij} = (1+|x|^{-1})^{4}\delta_{ij}  + \sigma_{ij},
\]
where $\partial_{I}\sigma_{ij} = O(|x|^{-2-|I|})$ for all multi-indices $I$ of length $|I|\leq 6$. We also assume that 
\[
x^{i}x^{j}\partial_i \partial_j R \geq 0 
\]
outside of a compact set. This condition integrates to yield 
\[
x^{i}\partial_{i}R \leq 0 \qquad\text{and} \qquad R \geq 0
\]
We now consider a sequence of connected closed stable constant mean curvature surfaces $\Sigma_{k}$ with
\[
r_{0}(\Sigma_{k})\to\infty, \qquad \area_g(\Sigma_k) \to \infty, \qquad\text{and} \qquad r_{0}(\Sigma_{k})H(\Sigma_{k})\to\infty.
\]
For $k$ large, we may find $\lambda>0$ and $\xi\in\RR^{3}$ large so that $\Sigma_{k} = \Sigma_{(\xi,\lambda)}$. By Theorem \ref{theo:far-off-center-LS}, 
\[
\frac{d}{ds}\Big|_{s=1} \area_{g}(\Sigma_{(s\xi,\lambda)}) = 0
\]
so that, by \eqref{lsradial}, 
\[
0 =  \frac{\pi}{105} \left( - 14 \, \lambda^5 |\xi| \partial_r \underline R - \lambda^7 |\xi| \partial_r \Delta \underline R + 144 \,   |\xi|^{-6}  \right) + O (\lambda^{-1} |\xi|^{-6}) + O (|\xi|^{-7}).
\]
It follows that 
\[
\partial_{r} \underline R = O(\lambda^{-5}|\xi|^{-7}) = o(\lambda^{-5}|\xi|^{-5}).
\]
Using this and \eqref{eq:R-concave-radial}, we may integrate in the in the radial direction to find that for $t\geq 0$,
\[
(\partial_{r}R) ( (1 + t)\lambda \, \xi) \geq \partial_{r} \underline R = o(\lambda^{-5}|\xi|^{-5}).
\]
Integrating this again, we find that
\[
\underline R \leq o(\lambda^{-4}|\xi|^{-4})t + R( (1 + t) \, \lambda \,  \xi ) \leq O(\lambda^{-4}|\xi|^{-4})(o(1)t + (1+t)^{-4}).
\]
Choosing $t$ judiciously we arrange for the term in parenthesis to be $o(1)$. We have proven that  
\[
\underline R = o(\lambda^{-4}|\xi|^{-4}).
\]
Now, considering the first variation of $\area_{g}(\Sigma_{(\xi,\lambda)})$ in directions orthogonal to $\xi$ as above. We obtain that the full derivative satisfies
\[
D\underline R = O(\lambda^{-5}|\xi|^{-7}).
\]
On the other hand, because $\partial_{r}\underline R = o(\lambda^{-5}|\xi|^{-5})$, Taylor's theorem combined with $\partial_{r}R \leq 0$ yields 
\[
\partial^{2}_{r}\underline R = o(\lambda^{-6}|\xi|^{-6}).
\]
Combining this with \eqref{eq:R-concave-radial}, we obtain  
\[
\partial^{3}_{r}\underline R = o(\lambda^{-7}|\xi|^{-7}).
\]
Similarly, combining the facts $R \geq 0$, $\underline R = o(\lambda^{-4}|\xi|^{-4})$, and $D\underline R = o(\lambda^{-5}|\xi|^{-5})$ with Taylor's theorem yields
\[
D^{2}\underline R \geq - o(\lambda^{-6}|\xi|^{-6}). 
\]
Similarly, we find that
\[
D^{2}\partial_{r}\underline R \leq o(\lambda^{-7}|\xi|^{-7}).
\]
Finally, since
\[
\partial_{r}\Delta \underline R = \Delta\partial_{r}\underline R - 2\, |\xi|^{-1}\lambda^{-1}\Delta \underline R + 2\, \lambda^{-1}|\xi|^{-1}\partial^{2}_{r}\underline R + 2\,  \lambda^{-2} |\xi|^{-2} \partial_{r}\underline R,
\]
we see that 
\[
\partial_{r}\Delta\underline R  \leq o(\lambda^{-7}|\xi|^{-7}).
\]
Returning to the radial first variation, we see that  
\[
0 \geq 14 \, \lambda^5 |\xi| \partial_r \underline R \geq  144 \, |\xi|^{-6} + \lambda^{-1} O(|\xi|^{-6}). 
\]
This contradiction completes the proof. 


\section{Proof of Theorem \ref{thm:homog-counterexample}}

As in the proof of Theorem \ref{thm:exist-outlying}, our strategy is parallel to the proof of Theorem 1 in \cite{Brendle-Eichmair:2014}, except that here we also exploit that the various terms in the reduced area functional $\xi \mapsto \area_g(\Sigma_{(\xi, \lambda)})$ have different orders in the regime where $\xi \to \infty$. \\

Consider $S: (0,\infty) \to (-\infty,0]$ a smooth function with 
\[
S^{(\ell)} = O(r^{-5-\ell})
\]
where $S^{(\ell)}$ is the $\ell$-th derivative. We define a smooth function $\varphi:(0,\infty)\to \RR$ by 
\[
\varphi(r) = \frac 1 r \int_{r}^{\infty}(\rho-r) \, \rho \, S(\rho) \, d\rho.
\]
Arguing as in Lemma \ref{lemm:order-falloff-varphi-counterexample}, we find that 
\[
\varphi^{(\ell)}(r) = O(r^{-5-\ell}).
\]
On the complement of a compact subset of $\RR^{3}$ we define a conformally flat Riemannian metric
\[
g = (1+1/r+\varphi(r))^{4}\bar g = (1+1/r)^{4}\bar g + O(1/r^{3}).
\]
Note that we can write
\[
g = (1+1/r)^{4}\bar g + T_{ij} + o(1/r^2)
\]
for $T_{ij} = 0$, so this is indeed of the form asserted in Theorem \ref{thm:homog-counterexample}. The scalar curvature satisfies
\[
R = - 8 \, (1+O(1/r))S(r).
\]
Now, fix $\chi \in C^{\infty}(\RR)$ with support in $[4,6]$ that is positive on $(4,6)$. Assume that $\chi'(5) = -1$. Define
\[
S(r) = - \sum_{k=0}^{\infty}10^{-5k} \chi(10^{-k}r).
\]
Note that $S^{(\ell)}(r) = O(r^{-5-\ell})$, as above. 

Consider $\xi \in \R^3$ with $|\xi| =  10^{k} \, t$ for $t\in [3,7]$. Then, taking $\lambda=10^{k}$, we have that
\begin{align*}
\area_{g}(\Sigma_{(10^{k},\xi)}) & = 4\pi \lambda^{2} - \frac{2\pi}{15} 10^{4k}R(\lambda^{20k}\xi) - \frac{\pi}{105} 10^{6k} (\Delta R) (\lambda^{20k}\xi) - \frac{8\pi}{35} |\xi|^{-6} + O(10^{-7k})\\
& = 4\pi \lambda^{2} + \frac{2\pi}{15} 10^{-5k}\chi(t) - \frac{8\pi}{35}10^{-6k} t^{-6} + O(10^{-7k}).
\end{align*}
Thus,
\[
\frac{d}{ds}\Big|_{s=1} \area_{g}(\Sigma_{(10^{k},s\xi)}) =  \frac{2\pi}{15} 10^{-6k}\chi'(t) + \frac{48\pi}{35} 10^{-7k} \, t^{-6} + O(10^{-8k}).
\]
For $t=7$, we have 
\[
\frac{d}{ds}\Big|_{s=1} \area_{g}(\Sigma_{(10^{k},s\xi)}) = \frac{48\pi}{35} 10^{-7k} \, 7^{-6} + O(10^{-8k}) > 10^{-5-7k}
\]
for sufficiently large $k$. On the other hand, for $t=5$, we have
\[
\frac{d}{ds}\Big|_{s=1} \area_{g}(\Sigma_{(10^{k},s\xi)}) =  -\frac{2\pi}{15} \, 10^{-6k} + \frac{48\pi}{35} 10^{-7k} \, 7^{-6} + O(10^{-8k}) <- 10^{-1-6k} .
\]
It follows that for some $t_{k} \in (5,7)$ and any $\xi_{k} \in \R^3$ with $|\xi_{k}| = 10^{k}\, t_{k}$, the surface  $\Sigma_{(10^{k},\xi_{k})}$ is a stable constant mean curvature  sphere. This completes the proof. 


\appendix


\section{Some integral expressions}\label{app:some-integral-expressions}


In this appendix, we recall several standard identities that are used in the proof of Theorem \ref{theo:far-off-center-LS}. 


\subsection{Integrals over \texorpdfstring{$B_{r} (0)$}{the ball}} 

Note that
\[
\int_{B_{r}(0)} (y^{i})^{2} = \frac 1 3 \int_{B_{r}(0)} |y|^{2} = \frac{4\pi}{15} r^{5} \qquad \text{ for all } i = 1, 2, 3.
\]
Thus, for a symmetric tensor $A_{ij}$ on $\R^3$, we have 
\[
\sum_{i, j } \int_{B_{r}(0)} A_{ij} y^{i}y^{j} =  \frac{4\pi}{15} r^{5} \, \tr  A.
\]
Similarly,
\[
\int_{B_{r}(0)} (y^{i})^{4} 
= \frac{4\pi}{35} r^{7}
\]
and for $i\not = j$,
\[
\int_{B_{r}(0)} (y^{i})^{2}(y^{j})^{2} 
= \frac{4\pi}{105} r^{7}
\]
For a totally symmetric tensor $B_{ijkl}$ on $\R^3$, we have that
\begin{align*}
\sum_{i,j,k,l}\int_{B_{r}(0)}B_{ijkl} y^{i}y^{j} y^{k}y^{l} & = \sum_{i} B_{iiii} \int_{B_{r}(0)} (y^{i})^{4} + 3\sum_{i\not =j} B_{iijj} \int_{B_{r}(0)} (y^{i})^{2}(y^{j})^{2} \\
& \qquad = \frac{4\pi}{35} r^{7} \Big(\sum_{i} B_{iiii} + \sum_{i\not =j} B_{iijj}  \Big)\\
& \qquad= \frac{4\pi}{35} r^{7} \sum_{i,j} B_{iijj} .
\end{align*}


\subsection{Integrals over \texorpdfstring{$S_{r}(0)$}{the sphere}}

Note that
\[
\int_{S_{r}(0)} (y^{i})^{2} = \frac{4\pi}{3} r^{4}. 
\]

It follows that, for a symmetric tensor $A_{ij}$ on $\R^3$, 
\[
\sum_{i, j} \int_{S_r(0)} A_{ij}y^{i}y^{j} = \frac{4\pi}{3} r^{4} \tr A.
\]

Similarly,
\begin{align*}
\int_{S_r(0)} (y^{i})^{4} &= \frac{4\pi}{5} r^{6} \qquad \text{ for all } i = 1, 2, 3,\\
\int_{S_r(0)} (y^{i})^{2}(y^{j})^{2} &= \frac{4\pi}{15} r^{6} \qquad \text{ for all } i \neq j.
\end{align*}

Thus, for a totally symmetric tensor $B_{ijkl}$ on $\R^3$, we have
\[
\int_{S_r(0)} B_{ijkl}y^{i}y^{j}y^{k}y^{l} = \frac{4\pi}{5} r^{6} B_{iijj}
\]

If $B_{ijkl}$ is symmetric in the first two slots and in the second two slots separately, we obtain 
\begin{align*}
\sum_{i, j, k, l} \int_{S_r(0)}B_{ijkl} y^{i}y^{j}y^{k}y^{l} & = \sum_{i} B_{iiii} \int_{S_r(0)} (y^{i})^{4} + \sum_{i}\sum_{j\not = i} B_{iijj} \int_{S_r(0)} (y^{i})^{2}(y^{j})^{2}\\
& \qquad +2 \sum_{i}\sum_{j\not = i} B_{ijij} \int_{S_r(0)} (y^{i})^{2}(y^{j})^{2}\\
& = \frac{4\pi}{15}r^{6} \Big( 3\sum_{i} B_{iiii}  +  \sum_{i}\sum_{j\not = i} B_{iijj} + 2 \sum_{i\not = j} B_{ijij}\Big) \\
& = \frac{4\pi}{15}r^{6} \Big( \sum_{i,j} B_{iijj} + 2 \sum_{i,j} B_{ijij}\Big).
\end{align*}

Finally, we have
\begin{align*}
\int_{S_r(0)} (y^{i})^{6} &= \frac{4\pi}{7} r^{8} \qquad \text{ for all } i = 1, 2, 3, \\
\int_{S_r(0)} (y^{i})^{4} (y^{j})^{2} &= \frac{4\pi}{35} r^{8} \qquad \text{ when } i \neq j, \\
\int_{S_r(0)} (y^{1})^{2}(y^{2})^{2}(y^{3})^{2} &= \frac{4\pi}{105} r^{8}.
\end{align*}

Assume now that the tensor $C_{ijklmn}$ on $\R^3$  is symmetric in the first four indices and, separately, in the last two indices. Then,
\begin{align*}
& \sum_{i,j,k,l,m,n} \int_{S_r(0)} C_{ijklmn} y^{i}y^{j}y^{k}y^{l}y^{m}y^{n}\\
& \qquad= \sum_{i} C_{iiiiii}  \int_{S_r(0)} (y^{i})^{6}\\
& \qquad\qquad + 6 \sum_{\substack{i, j \\ \text{distinct}}} C_{iijjjj} \int_{S_r(0)}(y^{i})^{2}(y^{j})^{4}\\
& \qquad\qquad+ \sum_{\substack{i, j \\ \text{distinct}}}  C_{iiiijj} \int_{S_r(0)}(y^{i})^{4}(y^{j})^{2}\\
& \qquad\qquad+ 3 \sum_{ \substack {i,j,k \\ \text{distinct}}} C_{iijjkk}  \int_{S_r(0)}(y^{i})^{2}(y^{j})^{2}(y^{k})^{2}\\
& \qquad\qquad+ 4 \sum_{\substack{i, j \\ \text{distinct}}}  C_{iiijij} \int_{S_r(0)}(y^{i})^{4}(y^{j})^{2}\\
& \qquad\qquad+ 12 \sum_{\substack{i, j, k \\ \text{distinct}}} C_{iijkjk}  \int_{S_r(0)}(y^{i})^{2}(y^{j})^{2}(y^{k})^{2}\\
& = \frac{4\pi}{35} r^{8}\Big( \sum_{i,j,k} C_{iijjkk} + 4 \sum_{i,j,k} C_{iijkjk} \Big).
\end{align*}


\subsection{Some useful integrals} \label{sec:someintegrals}
The following computations needed in the proof of Theorem \ref{theo:far-off-center-LS} are readily verified using the identities from the previous subsection.
\begin{align*}
\int_{S_r(0)} \Bigl( \frac{|a|^{2}|y|^{2} - 3\bangle{a,y}^{2}}{|a|^{5}}\Bigl)^{2} & = \frac{16\pi}{5} \frac{r^{6}}{|a|^{6}} \\
\int_{S_r(0)} \Bigl( \underline\sigma(y,y) - \frac 1 3 |y|^{2}\tr \underline \sigma \Bigl) \Bigl( \frac{|a|^{2}|y|^{2} - 3\bangle{a,y}^{2}}{|a|^{5}}\Bigl) & = 
\int_{S_r(0)}  \underline\sigma(y,y) \Bigl( \frac{|a|^{2}|y|^{2} - 3\bangle{a,y}^{2}}{|a|^{5}}\Bigl) \\
& = \frac{8\pi}{15} \frac{r^{6}}{|a|^{3}} \Big (  \tr\underline\sigma - 3 |a|^{-2}  \underline \sigma(a,a) \Big) \\
\int_{S_r(0)}\Big( \underline\sigma(y,y) - \frac 1 3 |y|^{2}\tr \underline \sigma \Big)^{2} 
& = \frac{8\pi}{45} r^{6} \Big (3  |\underline\sigma|^{2} - (\tr\underline\sigma)^{2} \Big) = \frac{8\pi}{15} r^{6}|\underline{\mathring\sigma}|^{2}.
\end{align*}

\bibliography{bib} 
\bibliographystyle{amsplain}
\end{document}